\documentclass[a4paper, 11pt]{article}

\usepackage{color}
\usepackage{graphicx}
\usepackage{amsmath}
\usepackage{amsfonts}
\usepackage{amssymb}
\usepackage{amsthm}
\usepackage{epstopdf}
\usepackage{algorithm}
\usepackage{setspace}
\usepackage[cm]{fullpage}
\usepackage{url}

\theoremstyle{remark}
\newtheorem{remark}{Remark}
\newtheorem{problem}{Problem}

\newcommand{\intO}{\int_\Omega}

\newcommand{\R}{\mathbb{R}}

\DeclareMathOperator*{\argmin}{argmin}

\title{Verification of functional a posteriori \\ error estimates for obstacle problem in 2D}
\author{P. HARASIM, J. VALDMAN
}

\begin{document}
\maketitle


\begin{abstract}
We verify functional a posteriori error estimate proposed by S. Repin \cite{NeRe} for a class of obstacle problems. The obstacle problem is formulated as a quadratic minimization problem with constrains equivalently formulated as a variational inequality. New benchmarks with known analytical solutions in 2D are constructed based on 1D benchmark introduced by P. Harasim and J. Valdman \cite{HaVa}. Numerical approximation of the obstacle problem is obtained by the finite element method using bilinear elements on a rectangular mesh. Error of the approximation is meassured in the energy norm and bounded from above by a functional majorant, whose value is minimized with respect to unknown gradient field discretized by Raviart-Thomas elements and Lagrange multipliers field discretized by piecewise constant functions. 
\end{abstract}

\section{Introduction}
Problems with obstacles often arise in continuum mechanics. Their mathematical models are formulated in terms of variational inequalities  \cite{GLT, HHNL}. Typically, numerical treatment of obstacle problems is obtained by the finite element method combined with methods developed for quadratic minimization problems with constrains.
It was traditionally tackled by the Uzawa method, the interior point method, the active set method with gradient splitting and the semi-smooth Newton method among others \cite{DO, Ul}. \\

A priori analysis providing asymptotic estimates of the quality of finite elements approximations converging toward the exact solution was studied for obstacle problems e.g. in \cite{BHR, FALK}. For the survey of the most important techniques in a posteriori analysis (such as residual, gradient averaging or equilibration methods) we refer to the monographs \cite{AiOd, BabSt, BanRa}. Particular a posteriori estimates for variational inequalities including a obstacle problem are reported e.g. in \cite{BHS, CaMe, ZVKG} among others. \\

Our goal is to verify guaranteed functional a posteriori estimates expressed in terms of functional majorants derived by Repin \cite{NeRe, ReGruyter}. The functional majorant upper bounds are essentially different with respect to known a posteriori error estimates mentioned above. The estimates are obtained with the help of variational (duality) method which was developed in \cite{Re1,Re2} for convex variational problems. The method was applied to various nonlinear models including those associated with variational inequalities \cite{Re4}, in particular problems with obstacles \cite{BuRe}, problems generated by plasticity theory \cite{FuRe, ReVa2} and problems with nonlinear boundary conditions \cite{ReVa}.\\

Three benchmarks with known analytical solution of the obstacle problem are considered in numerical experiments. For a known benchmark introduced in \cite{NSV} constructed on a square domain assuming non-zero Dirichlet boundary conditions and a constant obstacle, values of exact energy and Lagrange multipliers are added. Two new benchmarks in 2D are constructed on a ring domain assuming zero Dirichlet boundary conditions and either constant or spherical obstacle. The construction was inspired by a one dimensional benchmark from \cite{HaVa}. \\

Numerical testing is done by own Matlab code providing a bilinear approximation of the obstacle problems on uniform rectangles, Raviart-Thomas approximation of the gradient field and piecewise constant approximation of the Lagrange multiplier field. The code is vectorized in manner of \cite{RaVa} to provide a fast computation of finer mesh rectangulations. \\

Outline of the paper is as follows. In Section 2, we formulate a constrained minimization problem, a perturbed minimization problem and explain how to derive a functional a posteriori error estimate. 
Benchmarks with known analytical solution are discussed in Section 3. 
Numerical tests performed in Matlab are reported in Section 4. 
Additional details of the finite element implementation are mentioned in Appendix. 

\section{Obstacle problem and its functional a posteriori error estimate}
In the following, $\Omega\subset\mathbb{R}^2$
denotes a bounded domain with Lipschitz continuous boundary $\partial\Omega$, 
$V$ stands for the standard Sobolev space $H^{1}(\Omega)$ and $V_{0}$ denote its subspace $H^{1}_{0}(\Omega)$, consisting
of functions whose trace on $\partial\Omega$ is zero. We deal with the abstract obstacle problem, described by the following minimization problem:
\begin{problem}[Obstacle minimization problem] \label{prob1} Find $u \in K$ satisfying
\begin{equation}
J(u)=\inf_{v \in K} J(v),   \label{2rov1}
\end{equation}
where the energy functional reads
\begin{equation}
J(v) := \frac{1}{2}\intO\nabla v\cdot\nabla v \,\textrm{d}x-\intO fv\,\textrm{d}x      \label{2rov2}
\end{equation}
and the admissible set is defined as
$$
K:=\bigl\{v \in V_{0}:\,
v(x) \geq \phi(x) \mbox{ a.e.}\;\textrm{in}\; \Omega \},
$$
where a loading $f\in L^{2}(\Omega)$ and an obstacle $\phi\in V$, with $\phi(x)<0$  on $\partial \Omega$.
\end{problem}

Problem \ref{prob1} is a quadratic minimization problem with a convex constrain and the existence of its minimizer is guaranteed by the Lions-
Stampacchia Theorem \cite{LiSt}. It is equivalent to the following variational inequality: Find $u \in K$ such that
\begin{equation}
\intO \nabla u\cdot\nabla (v-u) \textrm{d}x\geq\intO f(v-u)\textrm{d}x \quad \mbox{for all }v \in K. \label{2varineq}
\end{equation}
The convex constrain $v\in K$ can be transformed into a linear term containing a new (Lagrange) variable in

\begin{problem}[Perturbed problem] \label{prob2} 
For given
$$
\mu \in \Lambda:=\left\{\mu\in L^2(\Omega): \mu\geq0\; \mbox{a.e.}\; \textrm{in}\; \Omega\right\}
$$
find $u_{\mu}\in V_{0}$ such that
\begin{equation} \label{2defpert}
J_{\mu}(u_{\mu})=\inf\limits_{v\in V_{0}} J_{\mu}(v),
\end{equation}
where the perturbed functional $J_{\mu}$ is defined as
\begin{equation}  \label{2Pert}
J_{\mu}(v):=J(v)-\intO \mu (v-\phi) \,\textrm{d}x.
\end{equation}
\end{problem}

\noindent Problems \ref{prob1} and \ref{prob2} are related and it obviously holds
\begin{equation} \label{2ineqJmuJ}
J_{\mu}(u_\mu) \leq J(u) \quad \mbox{for all }\mu \in \Lambda.
\end{equation}

\begin{remark}[Existence of optimal multiplier]
If $u$ has a higher regularity, 
\begin{equation}\label{smooth}
u\in V_{0}\cap H^{2}(\Omega),
\end{equation}
there exists an optimal multiplier
$\lambda \in \Lambda$ such that 
$u_{\lambda}=u$ and
$J_{\lambda}(u)= J(u)$. Moreover, it holds
\begin{equation} \label{podmlambda}
\lambda =-(\Delta u \,+f).
\end{equation}
For more details, see \cite[Lemma 2.1 and Remark 2.3]{HaVa}. 
\end{remark}
\begin{remark}[Higher regularity of solution]\label{rem:higher_regularity}
If $f$ and $\phi$ have higher regularity, e.g., if 
$f \in C(\overline \Omega)$ and $\phi \in C^1(\overline \Omega)$, 
then \eqref{smooth} holds and even $u \in C^1(\overline \Omega)$. For more details, see e.g. \cite{Gu, KS}.
\end{remark}

We are interested in analysis and numerical properties of the a posteriori error estimate of numerical
solution $v\in K$ to Problem \ref{prob1} in the energy norm
$$
\|v\|_{E}:=\left(\intO\nabla v\cdot\nabla v\,\textrm{d}x\right)^{\frac{1}{2}}.
$$
The following part is based on results of S. Repin et al. \cite{BuRe, NeRe, Re4}.
It is simple to see that
\begin{eqnarray}\label{energy_difference}
J(v)-J(u)=\frac{1}{2}\|v-u\|_{E}^2 +\intO \nabla u\cdot\nabla (v-u)\,\textrm{d}x-\intO f(v-u)\textrm{d}x 
\qquad \mbox{for all } v\in K
\end{eqnarray}
and \eqref{2varineq} implies an energy estimate
\begin{equation}\label{energyEstimate}
\frac{1}{2}\|v-u\|^{2}_{E}\leq J(v)-J(u) \quad \mbox{for all } v\in K.
\end{equation}
Quality of \eqref{energyEstimate} is tested in Section \ref{sec:numerics} for several problems with known exact solution $u\in K$ introduced in Section \ref{sec:benchmarks}. 
By using \eqref{2ineqJmuJ}, we obtain the estimate
\begin{equation}\label{eq9}
J(v)-J(u) \leq J(v)-J_{\mu}(u_{\mu})  \quad \mbox{for all }\mu \in \Lambda
\end{equation}
and S. Repin transformed \eqref{eq9} futher in a majorant estimate
\begin{equation}\label{majorantEstimate}
J(v)-J(u)\leq \mathcal{M}(v,f,\phi;\beta,\mu,\tau^{*}) \qquad \mbox{ for all } \beta>0, \mbox{ for all } \mu \in \Lambda, \mbox{ for all } \tau^{*}\in H(\Omega,{\rm{div}})
\end{equation}
in which the 'flux' variable $\tau^{*}$ is formulated in the space 
$$H(\Omega,{\rm{div}}):=\{\tau^{*}\in [L^2(\Omega)]^2:{\rm{div}}\,\tau^{*}\in L^{2}(\Omega)\}$$
well studied in various mixed problems. For more details on derivation of the majorant estimate, we refer to \cite{NeRe}. The right-hand side of \eqref{majorantEstimate} denotes a functional majorant
\begin{multline} \label{majorantForm}
\mathcal{M}(v,f,\phi;\beta,\mu, \tau^{*}) := \frac{1+\beta}{2} \|{\nabla v-\tau^{*}}\|^{2}_\Omega
+\frac{1}{2}\left(1+\frac{1}{\beta}\right) C_{\Omega}^2 \|{\rm{div}} \,\tau^{*}+f+\mu\|^{2}_\Omega
+\intO \mu(v-\phi)\textrm{d}x,
\end{multline}
where corresponding $L^2$ norms of vector and scalar arguments read
$$|| \tau^* ||_{\Omega} : =\left( \intO \tau^* \cdot \tau^* \,\textrm{d}x  \right)^{1/2} \qquad \mbox{for all } \tau^* \in L^2(\Omega,\mathbb{R}^{d}),$$ 
$$|| v ||_{\Omega} : =\left( \intO v^2 \,\textrm{d}x  \right)^{1/2} \qquad \mbox{for all } v \in L^2(\Omega)$$ 
and a constant $C_{\Omega}>0$ represents a smallest possible constant in the Friedrichs inequality $\| v \|_\Omega \leq C_{\Omega} \| \nabla v \|_\Omega$
valid for all $v\in V_{0}$. \\


\begin{remark}\label{rem}
If the assumption (\ref{smooth}) is fulfilled, there exist optimal majorant parameters $\tau_{\mathrm{opt}}^{*}=\nabla u$, $\mu_{\mathrm{opt}}=\lambda\in \Lambda$ and $\beta_{\mathrm{opt}}\rightarrow 0$ such that
the inequality in \eqref{majorantEstimate} changes to equality, i.e.
the majorant on right-hand side of (\ref{majorantEstimate}) defines the difference of energies\, $J(v)-J(u)$\, exactly
(see \cite[Lemma 3.4]{HaVa}). 
\end{remark}

\noindent Our goal is to find unknown optimal parameters $\beta_{\mathrm{opt}}$, $\mu_{\mathrm{opt}}$ and $\tau^{*}_{\mathrm{opt}}$ to keep the majorant estimate \eqref{majorantEstimate} as sharp as possible. 
\begin{problem}[Majorant minimization problem] \label{prob4}
Let $v\in K$, $f\in L^2(\Omega)$, $\phi<0$ be given. Find 
$$
(\beta_{\mathrm{opt}},\mu_{\mathrm{opt}},\tau^{*}_{\mathrm{opt}}) = \argmin \limits_{\beta,\mu,\tau^{*}}\mathcal{M}(v,f,\phi;\beta,\mu,\tau^{*}).
$$
over arguments $\beta > 0,  \mu \in \Lambda, \tau^{*}\in H(\Omega,{\rm{div}})$.
\end{problem}
\begin{algorithm}[ht]
\noindent Let $k=0$ and let $\beta_{k}>0$ and $\mu_{k}\in\Lambda$ be given. Then:
\begin{itemize}
\item[(i)] find $\tau^{*}_{k+1}\in H(\Omega,{\rm{div}})$ such that
$$
\tau^{*}_{k+1}=\argmin\limits_{\tau^{*}\in H(\Omega,{\rm{div}})}\mathcal{M}(v,f,\phi;\beta_k,\mu_k,\tau^{*}),
$$
\item[(ii)] find $\mu_{k+1}\in \Lambda$ such that
$$
\mu_{k+1}=\argmin\limits_{\mu\in \Lambda}\mathcal{M}(v,f,\phi;\beta_k,\mu,\tau^{*}_{k+1}),
$$
\item[(iii)] find $\beta_{k+1}>0$ such that
$$
\beta_{k+1}=\argmin\limits_{\beta > 0 }\mathcal{M}(v,f,\phi;\beta,\mu_{k+1},\tau^{*}_{k+1}),
$$
\item[(iv)] set $k:=k+1$ are repeat (i)-(iii) until convergence.
\item[(v)] output $\tau^{*}:=\tau^{*}_{k+1}$ and $\mu:=\mu_{k+1}$.
\end{itemize}
\caption{Majorant minimization algorithm}
\label{Alg:Majorant}
\end{algorithm}
We apply Algorithm \ref{Alg:Majorant} from \cite{HaVa} to solve Problem \ref{prob4}.
The minimization in step (i) is equivalent to the
following variational equation: Find $\tau^{*}_{k+1} \in H(\Omega,{\rm{div}})$ such that
\begin{multline}\label{minization1D_step_i}
(1+\beta_{k})\intO\tau^{*}_{k+1}\cdot w\textrm{d}x
+\left(1+\frac{1}{\beta_{k}}\right)\intO {\rm{div}} \,\tau^{*}_{k+1}{\rm{div}}\,w\textrm{d}x\\
=(1+\beta_{k})\intO \nabla v\cdot w\textrm{d}x
-\left(1+\frac{1}{\beta_{k}}\right)\intO(f+\mu_k){\rm{div}}\,w\textrm{d}x \quad \mbox{for all } w\in H(\Omega,{\rm{div}}).
\end{multline}
The minimization in step (ii) is equivalent to the variational inequality:
Find $\mu_{k+1}\in \Lambda$ such that
\begin{equation}  \label{minization1D_step_ii}
\intO\left[\left(1+\frac{1}{\beta_{k}}\right)\left({\rm{div}} \,\tau^{*}_{k+1}+\mu_{k+1}+f\right)+v-\phi\right](w-\mu_{k+1})\textrm{d}x\geq0
\quad \mbox{for all } w\in \Lambda.
\end{equation}
The minimization in step (iii) leads to the explicit relation
\begin{equation}  \label{minization1D_step_iii}
\beta_{k+1}=\frac{\|{\rm{div}} \,\tau^{*}_{k+1}+f+\mu_{k+1}\|_{L^{2}(\Omega)}}{ \| \nabla v-\tau^{*}_{k+1} \|_{[L^{2}(\Omega)]^2}} \,.
\end{equation}
Further details on implementation of Algorithm \ref{Alg:Majorant} are described in Subsection \ref{subsec:implementation_fem}.

\section{Benchmarks with known analytical solutions}\label{sec:benchmarks}
Three following benchmark problems provide exact solution of Problem \ref{prob1}. 

\subsection{Benchmark I: square domain, constant obstacle, nonzero Dirichlet BC}\label{subsec:benchamark_square} 
This benchmark is taken from \cite{NSV}.
Let us consider a square domain $\Omega=(-1,1)^2$ and prescribe a contact radius $R\in[0,1)$. For loading
\begin{equation*}
f(x,y)=\left\{
\begin{array}{lrl} -16(x^{2}+y^{2})+8R^2 & \quad\textrm{if} & \sqrt{x^2+y^2}>R\\
 -8(R^{4}+R^{2})+8R^{2}(x^{2}+y^{2}) & \quad\textrm{if} & \sqrt{x^2+y^2}\leq R, 
\end{array} \right.
\end{equation*}
it can be shown that 
\begin{equation*}
u(x,y)=\left\{
\begin{array}{lrll}
\left(\max\{x^{2}+y^{2}-R^{2},0\}\right)^2 & \quad \mbox{if} &(x,y) \in \Omega \\
\left(x^{2}+y^{2}-R^{2}\right)^2 & \quad \mbox{if} &(x,y) \in \partial\Omega.
\end{array} \right.
\end{equation*}
is the exact solution of Problem \ref{prob1} in case of the zero obstacle function $\phi=0$.
The corresponding energy reads
$$
J(u) = 192\left(\frac{12}{35}-\frac{28R^{2}}{45}+\frac{R^{4}}{3}\right)
-32R^{2}\left(\frac{28}{45}-\frac{4R^{2}}{3}+R^{4}\right)
+\frac{2}{3}\pi R^{8}.
$$ 
It is not difficult to show that
\begin{equation*}
\nabla u(x,y)=\left\{
\begin{array}{lrl}
4(x^{2}+y^{2}-R^{2})(x,y) & \quad\textrm{if} & \sqrt{x^2+y^2}>R\\
(0,0) & \quad\textrm{if} & \sqrt{x^2+y^2}\leq R. 
\end{array} \right.
\end{equation*}
With respect to (\ref{podmlambda}), the optimal multiplier reads
\begin{equation*}
\lambda(x,y)=\left\{
\begin{array}{lrl} 0 & \quad\textrm{if} & \sqrt{x^2+y^2}>R\\
8(R^{4}+R^{2})-8R^{2}(x^{2}+y^{2}) & \quad\textrm{if} & \sqrt{x^2+y^2}\leq R. 
\end{array} \right.
\end{equation*}

\subsection{Benchmarks defined on ring domain}\label{subsec:benchmarks_ring} 
We consider a ring domain
$\Omega:=\left\{(x,y)\in \mathbb{R}^{2}:x^{2}+y^{2}<1 \right\}$
and a constant negative loading $f$. In case of an inactive obstacle, Problem \ref{prob1} can be reduced to the following linear boundary value problem: Find a function
$u$
such that  
\begin{eqnarray}
-\Delta u=f&\rm{on}&\Omega  \label{motiveq}    \\
u=0&\rm{on}&\partial\Omega  \label{motivbound}.   
\end{eqnarray}
The solution of (\ref{motiveq})-(\ref{motivbound}) reads 
\begin{equation}\label{sol_exact_lin}
u(x,y) = \frac{f}{4}(1-x^2-y^2)   
\end{equation}
and\eqref{sol_exact_lin} minimizes the original energy functional (\ref{2rov2}) of Problem \ref{prob1} in the whole space $V_0$ (no obstacle constrain is respected)
and the corresponding energy reads
\begin{equation}\label{en_exact_lin}
J(u) = -\frac{\pi f^2}{16}.    
\end{equation}
In the case of an active and radially symmetric obstacle $\phi=\phi(r)$, where $r:=\sqrt{x^2+y^2}$ denotes a radius in polar coordinates, there will be an unknown contact domain radius $R\in (0,1)$ depending on the value of $f$ and the shape $\phi=\phi(r)$ only, see Figure \ref{fig:benchmark}. In Subsections \ref{subsub_constant} and \ref{subsub_spherical}, we consider cases of constant and spherical obstacles. Outside of the contact domain, where $r \in (R,1)$, the exact solution $u=u(r)$ will be radially symmetric and satisfies the equation (\ref{motiveq}) transformed to polar coordinates and modified boundary conditions
\begin{eqnarray}
\frac{\partial^2 u}{{\partial r}^2}+\frac{1}{r}\frac{\partial u}{\partial r}&=&-f \qquad \mbox{for } r \in (R,1)   \label{TransLapl} \\
\vspace{0.1cm} 
u(R)&=&\phi(R) \\
u(1)&=&0.  \label{cond}     
\end{eqnarray}
The solution of (\ref{TransLapl})-(\ref{cond}) reads  
\begin{equation}
u(r)=\frac{4\phi(R) +fR^2-f}{4\ln R}\ln r-\frac{fr^2}{4}+\frac{f}{4}.  \label{sol1D}
\end{equation}

\subsubsection{Benchmark II: ring domain, constant obstacle and zero Dirichlet BC}\label{subsub_constant}
This benchmark generalizes 1D benchmark from \cite{HaVa} into 2D.
We consider a constant negative obstacle function $\phi$.
It follows from (\ref{sol_exact_lin}), the obstacle is inactive for smaller loadings satisfying
$|f| < 4|\phi|$. Thus, the exact solution of Problem \ref{prob1} is given by (\ref{sol_exact_lin}) and
the corresponding energy by (\ref{en_exact_lin}).
The obstacle is active if $|f| \geq 4|\phi|$ and, with respect of (\ref{sol1D}), the solution of Problem \ref{prob1} reads
\begin{equation*}
u(x,y)=\left\{
\begin{array}{lrl}
\frac{f}{4}\left(1-x^2-y^2\right) + A_{c}\ln\sqrt{x^2+y^2} & \quad\textrm{if} & \sqrt{x^2+y^2}>R\\
\phi & \quad\textrm{if} & \sqrt{x^2+y^2}\leq R, 
\end{array} \right.
\end{equation*}
where
$$
A_{c}=\frac{4\phi+fR^2-f}{4\ln R}.
$$
An unknown contact radius $R \in (0, 1) $ follows from the condition of continuity of the first derivative $\frac{\partial u(r)}{\partial r} \big|_{r=R}=0$, i.e., from the nonlinear equation
\begin{equation}\label{condR_c}
R^{2}(1-2\ln R)=1-\frac{4\phi}{f}.  
\end{equation}
The corresponding energy can be expressed as
$$
J(u) = \frac{\pi f^{2}R^4}{16}-\frac{\pi (\phi-\frac{f}{4})^2}{16\ln R}-\frac{\pi fR^{2}(\phi-\frac{f}{4})}{2\ln R}-\frac{\pi f^{2}R^{4}}{16\ln R}-\frac{\pi f^2}{16}.    
$$
It is not difficult to show that
\begin{equation*}
\nabla u(x,y)=\left\{
\begin{array}{lrl}
\left(\frac{A_{c}x}{x^2+y^2}-\frac{fx}{2},\frac{A_{c}y}{x^2+y^2}
-\frac{fy}{2}\right) & \quad\textrm{if} & \sqrt{x^2+y^2}>R\\
(0,0) & \quad\textrm{if} & \sqrt{x^2+y^2}\leq R. 
\end{array} \right.
\end{equation*}
With respect to (\ref{podmlambda}), the optimal multiplier reads
\begin{equation*}
\lambda(x,y)=\left\{
\begin{array}{crl}
0 & \quad\textrm{if} & \sqrt{x^2+y^2}>R\\
-f & \quad\textrm{if} & \sqrt{x^2+y^2}\leq R, 
\end{array} \right.
\end{equation*} 
so that it is a piecewise constant function.

\subsubsection{Benchmark III: ring domain, spherical obstacle and zero Dirichlet BC} \label{subsub_spherical}
\begin{figure}
\setlength{\unitlength}{3pt}
\begin{picture}(45,100)(-50,-40)
\put(-10,30){\vector(1,0){100}}
\put(33,36){\vector(0,1){15}}
\put(33,-35){\line(0,1){65}}
\put(0,35){\vector(0,-1){5}}
\put(3,35){\vector(0,-1){5}}
\put(6,35){\vector(0,-1){5}}
\put(9,35){\vector(0,-1){5}}
\put(12,35){\vector(0,-1){5}}
\put(15,35){\vector(0,-1){5}}
\put(18,35){\vector(0,-1){5}}
\put(21,35){\vector(0,-1){5}}
\put(24,35){\vector(0,-1){5}}
\put(27,35){\vector(0,-1){5}}
\put(30,35){\vector(0,-1){5}}
\put(33,35){\vector(0,-1){5}}
\put(36,35){\vector(0,-1){5}}
\put(39,35){\vector(0,-1){5}}
\put(42,35){\vector(0,-1){5}}
\put(45,35){\vector(0,-1){5}}
\put(48,35){\vector(0,-1){5}}
\put(51,35){\vector(0,-1){5}}
\put(54,35){\vector(0,-1){5}}
\put(57,35){\vector(0,-1){5}}
\put(60,35){\vector(0,-1){5}}
\put(63,35){\vector(0,-1){5}}
\put(66,35){\vector(0,-1){5}}

\multiput(49,18)(0,2){4}{\line(0,1){1}}

\put(89,27){r}
\put(66,27){1}
\put(-3,27){-1}
\put(12,36){$f=\mathrm{const}$}
\put(30,49){z}
\put(61,19){u(r)}
\put(33,21){$\phi_{\mathrm{max}}$}
\put(33,-31){$\phi_{\mathrm{max}}-\rho$}
\put(41,-7){$\rho$}
\put(56,9){$\phi(r)$}

\qbezier(33,-11)(36,-11)(39,-12)
\put(34,-17){$\psi$}

\put(49,31){\line(0,-1){2}}
\put(48,26){R}

\qbezier(33,20)(52,20)(66,10)
\qbezier(0,10)(14,20)(33,20)

\qbezier(0,30)(6,16)(17,18)
\qbezier(49,18)(60,16)(66,30)

\put(33,27){0}
\put(33,-30){\vector(1,3){16}}
\end{picture}
\caption{Benchmark setup: Constant forces $f$ pressing continuum against a lower spherical obstacle $\phi$ and the exact displacement $u$. The contact radius $R$ is determined by the formula  $R=\rho\sin\psi$ for some angular parameter $\psi$ and a prescribed value of the sphere radius $\rho$.}
\label{fig:benchmark}
\end{figure}
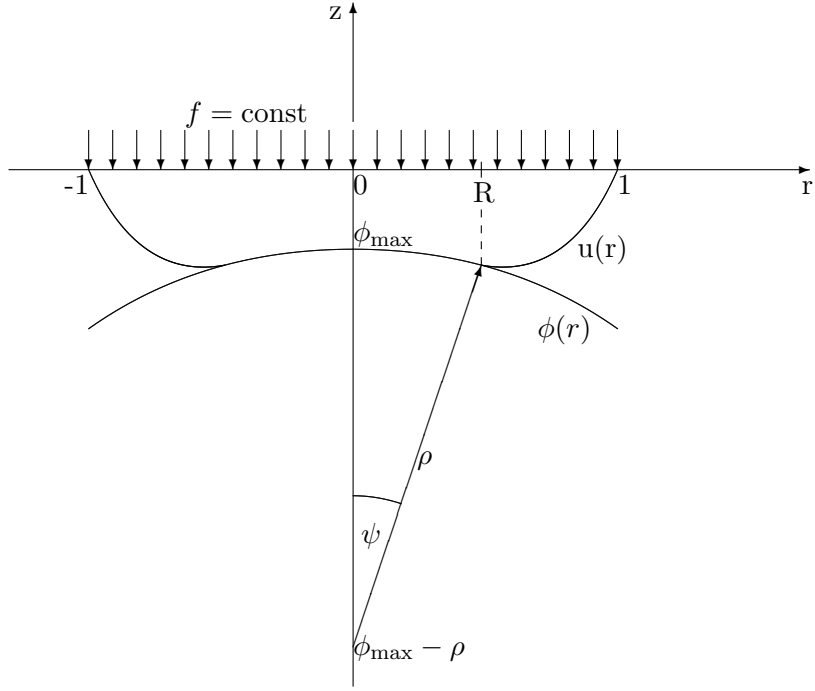
In this benchmark, we replace the constant obstacle $\phi$ 
by a sphere of the radius $\rho\geq1$. An obstacle function is defined as
$$
\phi(x,y)=\phi_{\mathrm{max}}-\rho+\sqrt{\rho^2-x^2-y^2}, \qquad (x,y)\in \Omega
$$ 
where $\phi_{\mathrm{max}}$ is a prescribed negative constant meaning the maximal value of the obstacle function. 
Analogically to the previous bechmark with the constant obstacle, for smaller loadings satisfying $|f| < 4|\phi_{\mathrm{max}}|$
we have a linear problem with inactive obstacle, the solution of Problem \ref{prob1} is given by (\ref{sol_exact_lin}) and corresponding energy
by (\ref{en_exact_lin}). The obstacle is active if 
$|f| \geq 4|\phi_{\mathrm{max}}|$ holds.
Then, the solution of Problem \ref{prob1} reads
\begin{equation*}
u(x,y)=\left\{
\begin{array}{lrl}
\frac{f}{4}\left(1-x^2-y^2\right) + A_{s}\ln\sqrt{x^2+y^2} & \quad\textrm{if} & \sqrt{x^2+y^2}>R\\
\phi_{\mathrm{max}}-\rho+\sqrt{\rho^2-x^2-y^2} & \quad\textrm{if} & \sqrt{x^2+y^2}\leq R, 
\end{array} \right. 
\end{equation*}
where 
$$
A_{s}=\frac{4(\phi_{\mathrm{max}}-\rho+\rho\cos\psi)+fR^{2}-f}{4\ln R}.
$$
An unknown contact radius $R=\rho\sin\psi$ for some angular parameter $\psi \in (0, \arcsin \frac{1}{\rho})$ (see Figure \ref{fig:benchmark} for details) follows from the condition of continuity of the first derivative
$
\frac{\partial u(r)}{\partial r} \big|_{r=R}= -\tan\psi,
$
i.e., from the solution of the nonlinear equation 
\begin{equation}\label{condR_s}
\frac{4(\phi_{\mathrm{max}}-\rho +\rho\cos\psi)+f\rho^2\sin^2\psi-f}{4\rho\sin\psi\ln(\rho\sin\psi)}
-\frac{f\rho\sin\psi}{2}=-\tan\psi.  
\end{equation}
The corresponding energy \eqref{2rov2} can be decomposed as  
$$J(u)=J_{1}(u)+J_{2}(u), $$
where the first term related to the contact domain 
$\Omega_{\ominus}^{u}=\left\{(x,y)\in \mathbb{R}^{2}:x^{2}+y^{2}\leq R \right\}$ reads
$$
J_{1}(u)= \frac{1}{2}\int_{\Omega_{\ominus}^{u}}\frac{x^{2}+y^{2}}{\rho^2-x^2-y^2}\,\textrm{d}x
-\int_{\Omega_{\ominus}^{u}}f\left(\phi-\rho+\sqrt{\rho^2-x^2-y^2}\right)\,\textrm{d}x
$$
and the second term related to the remaining part $\Omega_{0}^{u}:=\Omega \,\backslash \,\Omega_{\ominus}^{u}$ reads 
$$
J_{2}(u)= \frac{1}{2}\int_{\Omega_{0}^{u}}\left[\frac{A_{s}^{2}}{x^2+y^2}-A_{s}f+\frac{f^{2}}{4}\left(x^2+y^2\right)\right]\,\textrm{d}x
-\int_{\Omega_{0}^{u}}\left[\frac{f^{2}}{4}\left(1-x^2-y^2\right)+A_{s}f\ln\sqrt{x^2+y^2}\right]\,\textrm{d}x.
$$
They can be expressed as
$$
J_{1}(u)=-\frac{\pi\rho^2}{2}[\sin^{2}\psi+\ln(\cos^{2}\psi)]
-\pi f\rho^2(\phi-\rho)\sin^{2}\psi-\frac{2\pi f\rho^3}{3}\left(1-\sqrt{(1-\sin^{2}\psi)^3}\right)
$$
and 
$$
J_{2}(u)=-\pi A_{s}^{2}\ln R-\frac{\pi f(2A_{s}-1)(1-R^{2})}{4}-\frac{A_{s}\pi fR^{2}(1-2\ln R)}{2}
+\frac{3\pi f^{2}(1-R^{4})}{16}.
$$
It is not difficult to show that
\begin{equation*}
\nabla u(x,y)=\left\{
\begin{array}{lrl}
\Bigl(\frac{A_{s}x}{x^2+y^2}-\frac{fx}{2},
\frac{A_{s}y}{x^2+y^2}-\frac{fy}{2}\Bigr) & \textrm{if} & \sqrt{x^2+y^2}>R\\
\Bigl(-\frac{x}{\sqrt{\rho^2-x^2-y^2}},-\frac{y}{\sqrt{\rho^2-x^2-y^2}}\Bigr) & \textrm{if} & \sqrt{x^2+y^2}\leq R. 
\end{array} \right.
\end{equation*}
With respect to (\ref{podmlambda}), the optimal multiplier reads
\begin{equation*}
\lambda(x,y)=\left\{
\begin{array}{crl}
0 & \quad\textrm{if} & \sqrt{x^2+y^2}>R\\
\frac{2\rho^{2}-x^2-y^2}{(\rho^2-x^2-y^2)^{\frac{3}{2}}}-f & \quad\textrm{if} & \sqrt{x^2+y^2}\leq R. 
\end{array} \right.
\end{equation*} 
\begin{remark}
Existence of first continuous derivatives $
\frac{\partial u(r)}{\partial r} \big|_{r=R}
$ resulting in conditions \eqref{condR_c} and \eqref{condR_s} to determine the contact radius $R$ can be justified by Remark \ref{rem:higher_regularity}. 
\end{remark}

\section{Numerical experiments}\label{sec:numerics}
We verify the energy estimate \eqref{energyEstimate} and the majorant estimate \eqref{majorantEstimate} on three introduced benchmarks. 
Numerical experiments are based on an own implementation of the finite element method on uniform rectangular meshes in two dimensions: bilinear nodal basis functions for the solution $v$ of the obstacle problem (Problem \ref{prob1}) and Raviart-Thomas basis functions for construction the flux $\tau$ in the majorant minimization problem (Problem \ref{prob4}). Lagrange multiplies $\mu$ are constructed as piecewice constant funcions. A MATLAB code is available for download as a package {\it Obstacle problem
in 2D and its a posteriori error estimate} at Matlab Central under {\small\url{http://www.mathworks.com/matlabcentral/fileexchange/authors/37756}}. The implementation is based on vectorization techniques of \cite{RaVa} and works fast for finer rectangular meshes.

\subsection{Approximations by the finite element method}\label{subsec:implementation_fem}
Discrete approximation $v$ of the solution $u$ of Problem \ref{prob1} is expressed by a linear combination
$$v =\sum_{j=1}^{N}v_{j}\psi_{j}$$ of nodal bilinear functions $\psi_j$, where $N$ denotes a number of nodes of a considered rectangulation $\mathcal{T}$. Nodal components are assembled in a nodal (collumn) vector $${\bf v}=(v_1, \dots, v_N)^T.$$
Let us assume that $N-N_D$ internal nodes are ordered first and $N_D$ boundary Dirichlet nodes last. Then, ${\bf v}$ solves a quadratic minimization problem
\begin{equation}\label{prog:quadprog}
{\bf v}=\argmin_{w_i \geq \varphi_i, w_j = u_j} \left\{ \frac{1}{2} {\bf w} ^T {\bf K}^{BIL} {\bf w} - {\bf b}^T {\bf w} \right\},
\end{equation}
for $i \in \{1,\dots, N-N_D \}$ and $j \in \{N-N_D+1,\dots, N \}$. Here $\varphi_i$ and $u_j$ denote nodal obstacle values and nodal Dirichlet boundary values. A stifness matrix 
${\bf K}^{BIL}\in \R^{N \times N}$ and a discretized loading (collumn) vector ${\bf b} = (b_1, \dots, b_N)^T \in \R^N$ are constructed elemementwise as
\begin{equation}\label{eq:b_definition}
 ({\bf K}^{BIL})_{ij}=\int_{\Omega} \nabla  \psi_i \cdot \nabla  \psi_j \,\textrm{d}x, \qquad b_{i}=\int_{\Omega} f \psi_i \,\textrm{d}x 
\end{equation}
for $i, j \in \{1,\dots, N \}$. In all quadratures related to $f$ function, $f$ is replaced by a piecewise constant function $\overline{f}$ computed as the average of four nodal values on a rectangle. We apply the built-in Matlab function {\it quadprog} to solve \eqref{prog:quadprog}. 

A discretized version of Algorithm \ref{Alg:Majorant} is applied for the minimization of the functional majorant. 
The minimal argument $\tau^{*}_{k+1} \in H(\Omega,{\rm{div}})$ in step (i) of Algorithm \ref{Alg:Majorant} is expressed by a linear combination
$$\tau^{*}_{k+1}=\sum_{j=1}^{M}y_{j}\eta_{j},$$ 
of edge Raviart Thomas vector functions $\eta_{j}$, where $M$ denotes a number of rectangulation edges. Then, a coefficient (collumn) vector 
$${\bf y}=(y_1,\dots,y_M)^T$$ solves (see \eqref{minization1D_step_i}) a linear system of equations
\begin{equation}\label{eq:system_y}
\left[(1+\beta_k){\bf M}^{RT0}+\left(1+\frac{1}{\beta_k}\right){\bf K}^{RT0}\right]{\bf y}=(1+\beta_k){\bf c}-\left(1+\frac{1}{\beta_k}\right){\bf d}.
\end{equation}
Here, a stifness matrix ${\bf K}^{RT0} \in \R^{M \times M}$ and a mass matrix ${\bf M}^{RT0} \in \R^{M \times M}$ are constructed as
$$({\bf K}^{RT0})_{ij}=\int_{\Omega} {\rm{div}} \eta_i  \; {\rm{div}} \eta_j \,\textrm{d}x,   
\qquad ({\bf M}^{RT0})_{ij}=\int_{\Omega}  \eta_i \cdot \eta_j \,\textrm{d}x
$$ 
for $i, j \in \{1,\dots, M \}$
and ${\bf c} =(c_1, \dots, c_M)^T \in \R^M$ and ${\bf d} =(d_1, \dots, d_M)^T\in \R^M$ are (collumn) vectors constructed as
$$
c_i=\intO\nabla v\cdot\eta_{i}\textrm{d}x, \quad d_i=\intO(f+\mu_{k}){\rm{div}}\eta_{i}\,\textrm{d}x
$$
for $i \in \{1,\dots, M \}$. No boundary conditions are imposed on ${\bf y}$ in \eqref{eq:system_y} since the discrete solution $v$ satisfies Dirichlet boundary conditions only. The minimal argument $\mu_{k+1} \in \Lambda_{h}$ in step (ii) of Algorithm \ref{Alg:Majorant} is searched in the finite dimensional space $\Lambda_{h} \subset \Lambda$ of piecewise constant functions 
. It is computed locally on every rectangle from the formula
\begin{equation}
\mu_{k+1} 
=\left[-{\rm{div}}\,\tau^{*}_{k+1}-\overline{f}
-\frac{\overline{v} -\overline{\phi}}
{C_{\Omega}^{2}\left(1+\frac{1}{\beta_{k}}\right)}\right]^+,
\end{equation}
where $\overline{v}, \overline{\phi}$ represent averaged rectangular values computed as the average of four nodal values on a rectangle and $[\cdot]^+=\max\{ 0,\cdot \}$ denotes the maximum operator. \\

\begin{remark}
In all numerical experiments, only two iterations of Algorithm \ref{Alg:Majorant} were applied in which we set $\beta=1$ and $\mu_0$ is provided from the quadratic programming function {\it quadprog}. Without a good initial iteration $\mu_0$, the number of iterations would be significantly higher as demonstrated in \cite{HaVa}.
\end{remark}

\begin{remark}
Exact forms of local finite element matrices are reported in Appendix (Section \ref{sec:appendix}).
\end{remark}

\subsection{Computational results.}
We verify the energy estimate \eqref{energyEstimate} and the majorant estimate \eqref{majorantEstimate} for all three introduced benchmarks discretized on uniform rectangular meshes $\mathcal{T}_h$ generated by the uniform mesh parameter
$$h \in \left\{ \frac{1}{2},\frac{1}{4},\frac{1}{8},\frac{1}{16}, \frac{1}{32}, \frac{1}{64} \right\}.$$
The error of the discrete approximation $v_h$ is approximately evaluated by a quadratic form
\begin{equation}\label{eq:error_estimate}
\|v_h-u\|^{2}_{E} \approx ({\bf v}_h - {\bf u}_h)^T {\bf K}^{BIL}_h ({\bf v}_h - {\bf u}_h)
\end{equation}
using discrete solution ${\bf v}_h$ and nodal intepolation ${\bf u}_h$ of the exact solution $u$ and a stifness matrix matrix ${\bf K}^{BIL}_h$ on a considered mesh $\mathcal{T}_h$. The error value is further improved by evaluations of quadratic forms
\begin{equation}\label{eq:error_estimate2}
\|v_h-u\|^{2}_{E} \approx (\mathcal P_{h/2}({\bf v}_{h}) - {\bf u}_{h/2})^T {\bf K}^{BIL}_{h/2} (\mathcal P_{h/2}({\bf v}_{h}) - {\bf u}_{h/2})
\end{equation}
and 
\begin{equation}\label{eq:error_estimate3}
\|v_h-u\|^{2}_{E} \approx (\mathcal P_{h/4}({\bf v}_{h}) - {\bf u}_{h/4})^T {\bf K}^{BIL}_{h/4} (\mathcal P_{h/4}({\bf v}_{h}) - {\bf u}_{h/4})
\end{equation}
using nodal prolongation operators $\mathcal P_{h/2}$ and $\mathcal P_{h/4}$ to once and even twice uniformly refined rectangular meshes $\mathcal{T}_{h/2}$ and $\mathcal{T}_{h/4}$ and nodal interpolations ${\bf u}_{h/2}$ and ${\bf u}_{h/4}$.
Postprocessing of \eqref{eq:error_estimate2}, \eqref{eq:error_estimate3} requires extra memory resources but it proves important to keep a proper inequality sign in the energy estimate \eqref{energyEstimate}. The value of \eqref{eq:error_estimate3} is used in all convergence figures and all values \eqref{eq:error_estimate}, \eqref{eq:error_estimate2}, \eqref{eq:error_estimate3} are prompted in the code run for illustration. 

\begin{figure}
\center
\begin{minipage}{0.49\textwidth}
\center
\includegraphics[width=0.8\textwidth]{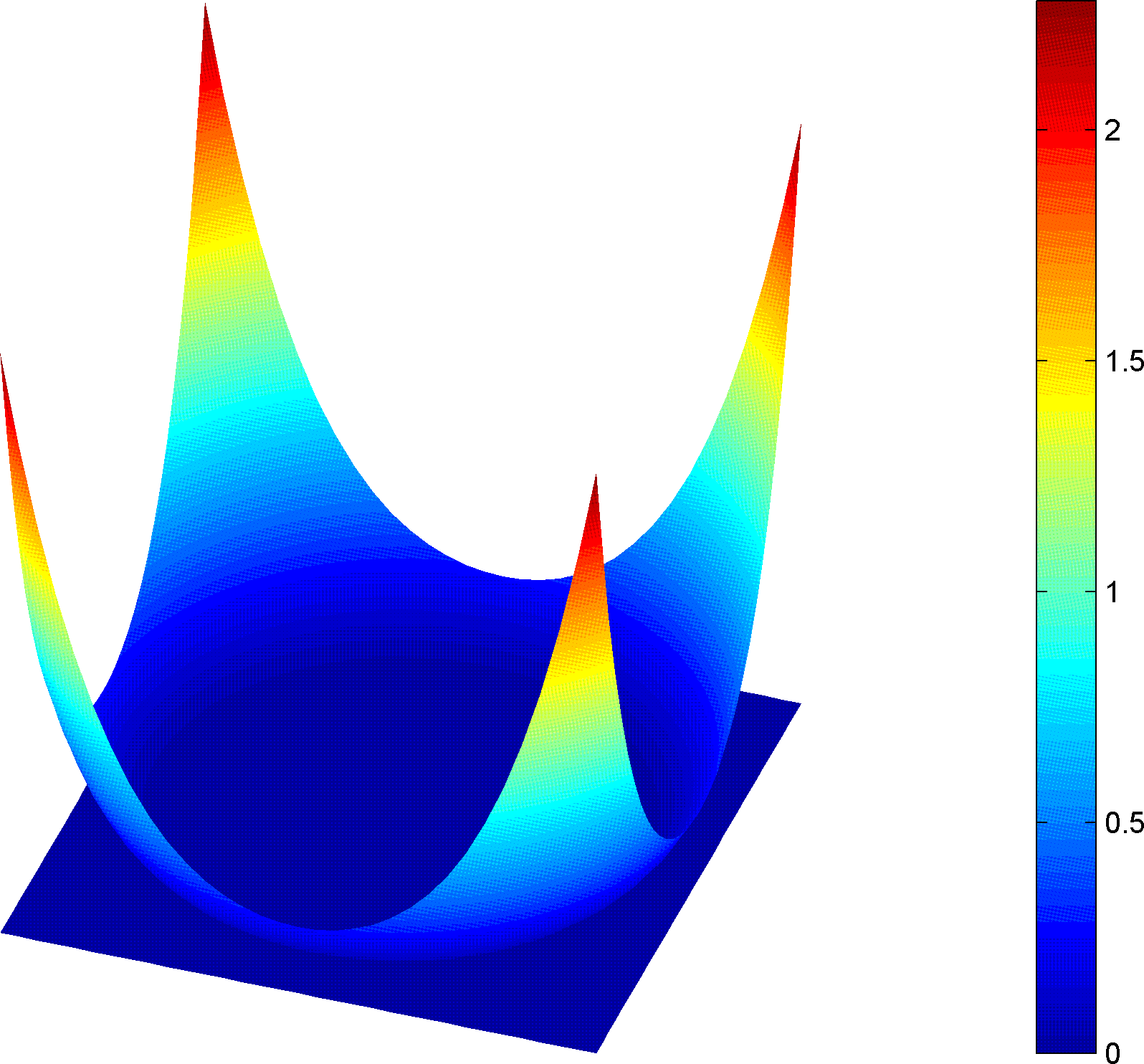}
\end{minipage}
\begin{minipage}{0.49\textwidth}
\center
\includegraphics[width=0.8\textwidth]{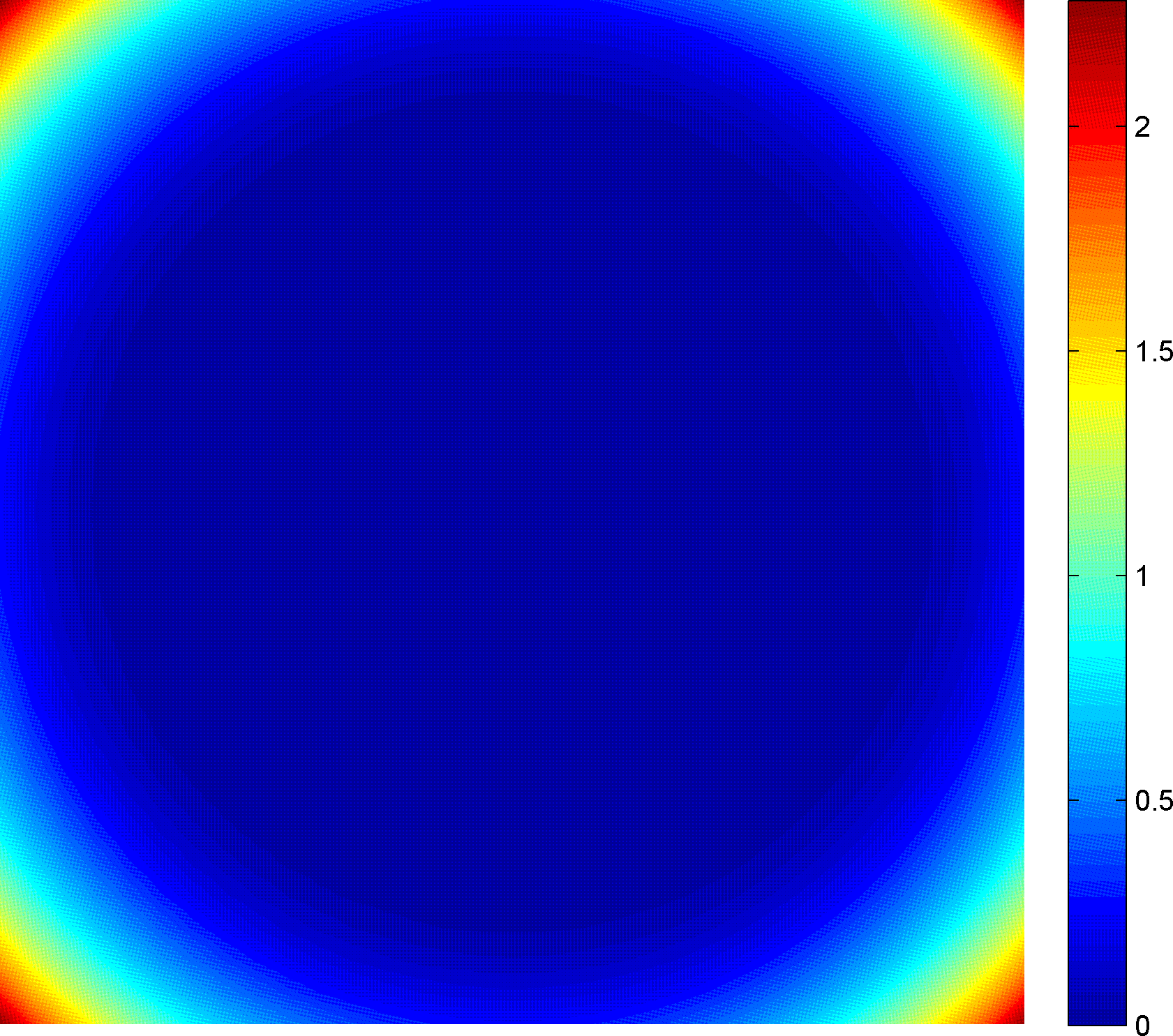}
\end{minipage}
\caption{Discrete solution $v$ of the obstacle problem and the lower obstacle $\phi$ (left) and its rotated view (right). The dark blue color indicates the contact domain.}
\label{benchmark_square_discrete_solution}
\begin{minipage}{0.49\textwidth}
\center
\includegraphics[width=0.8\textwidth]{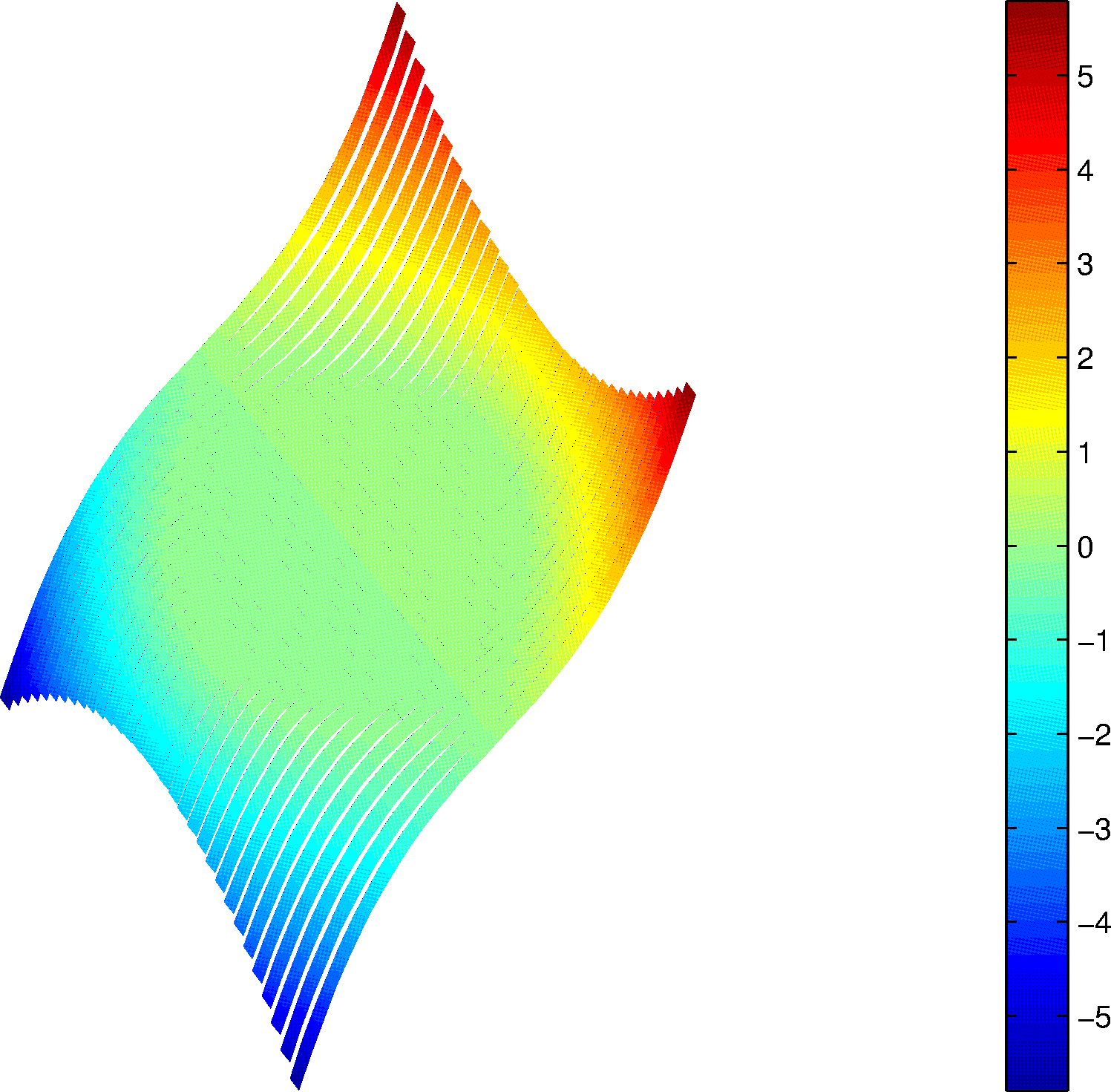}
\vspace{0.2cm} \\
\includegraphics[width=0.8\textwidth]{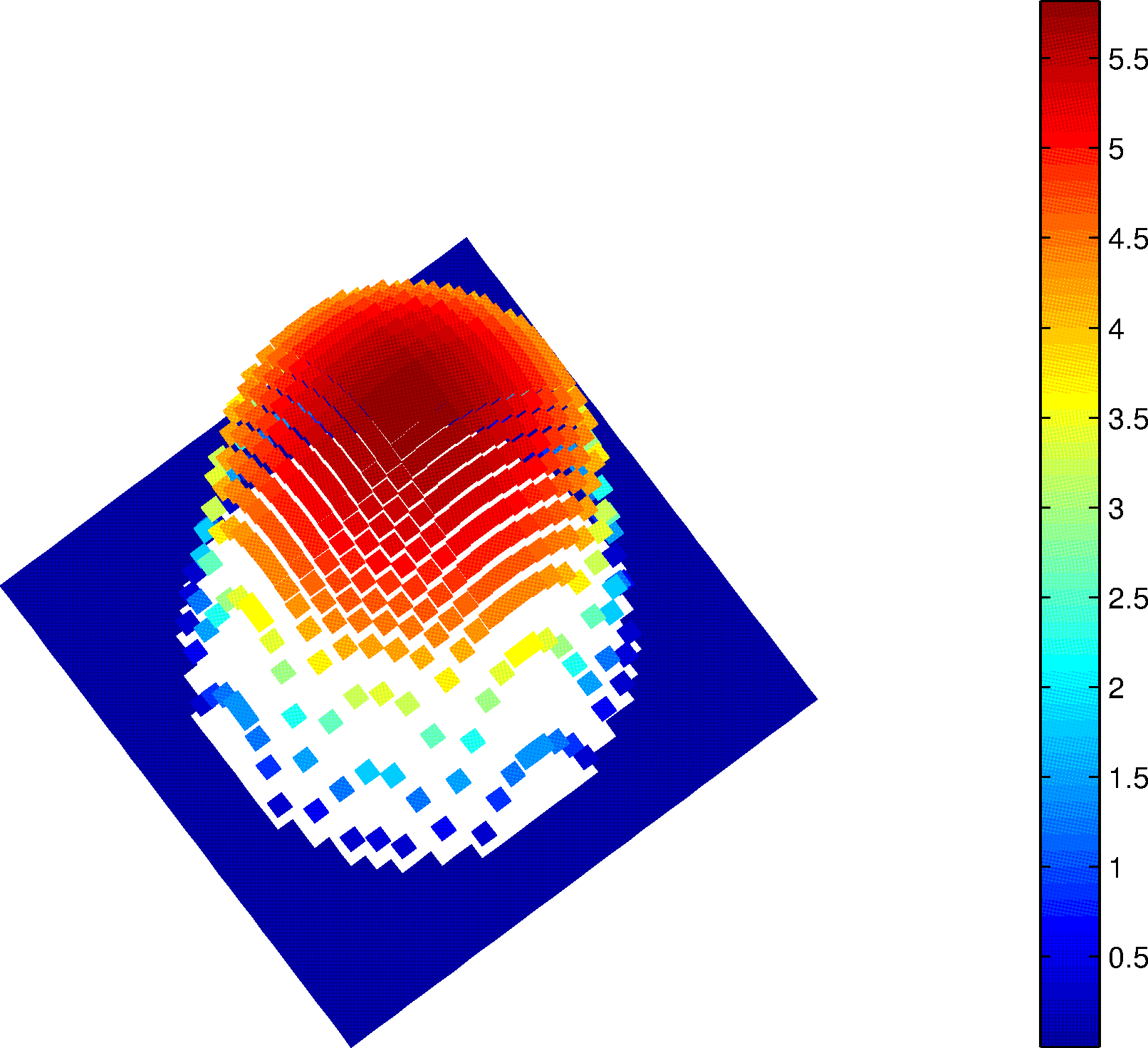}
\end{minipage}
\begin{minipage}{0.49\textwidth}
\center
\includegraphics[width=0.8\textwidth]{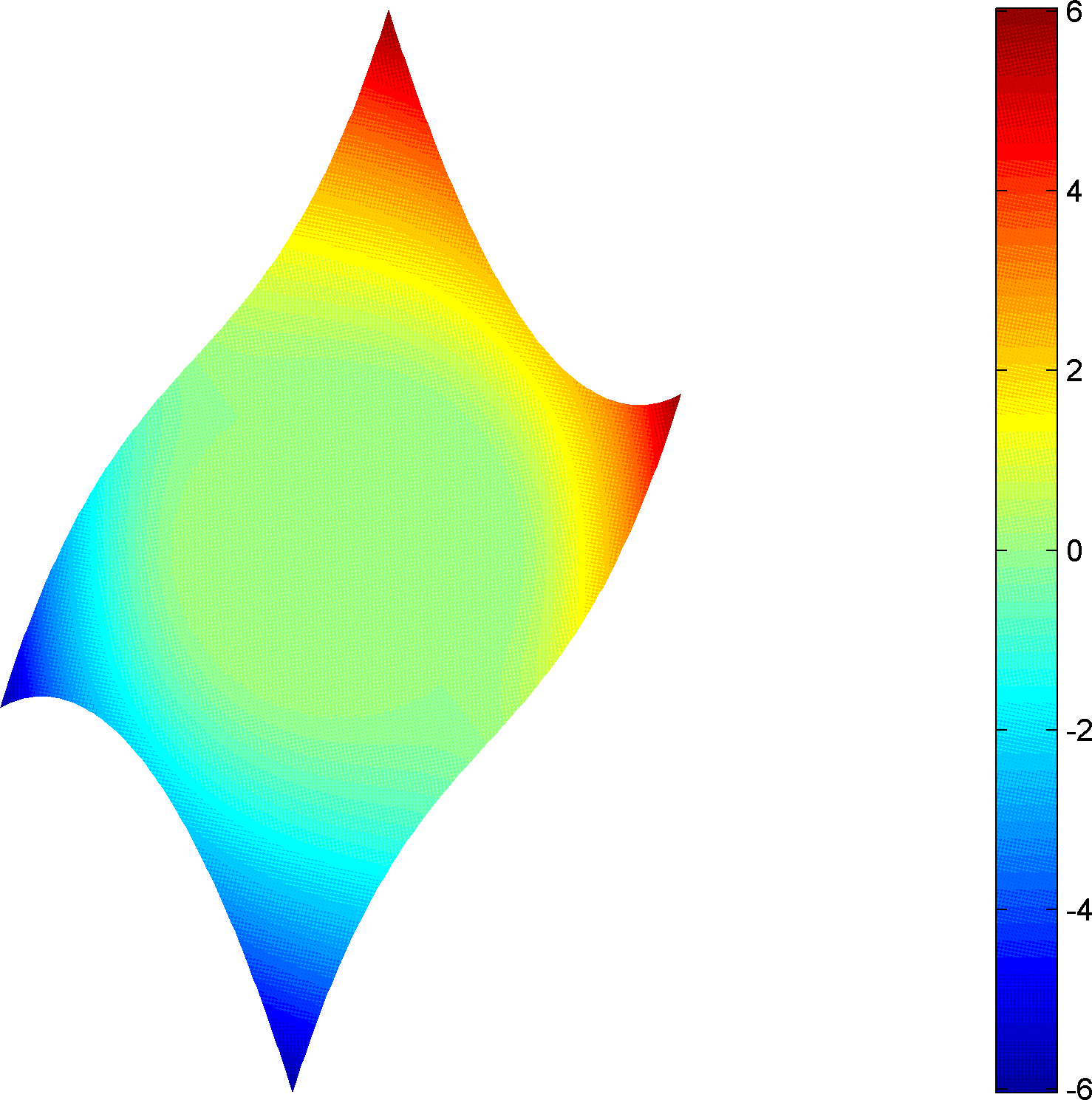}
\vspace{0.2cm} \\
\includegraphics[width=0.8\textwidth]{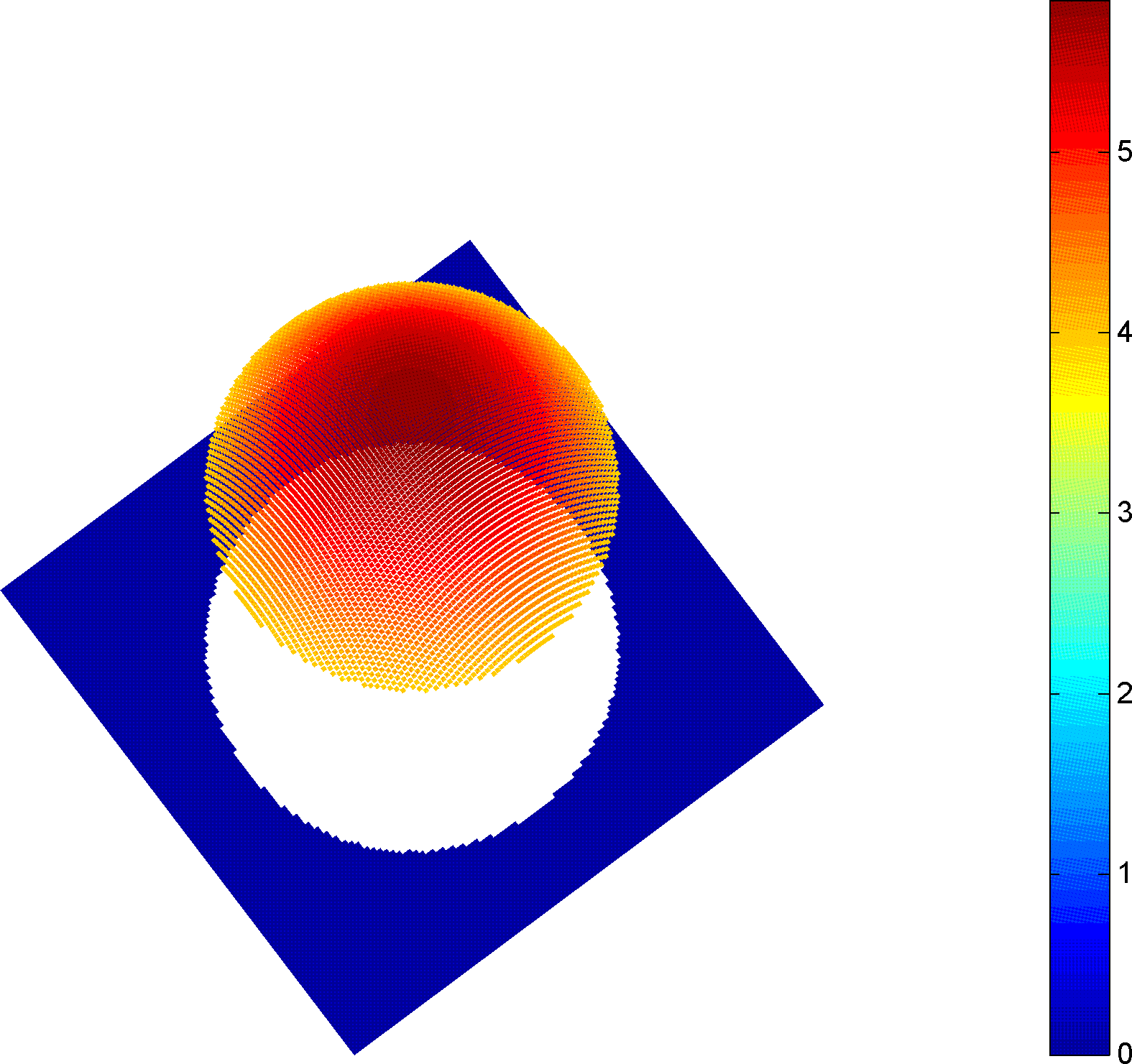}
\end{minipage}
\caption{Discrete flux x-component $\tau_x^*$ (top left) and discrete multiplier $\mu$ (bottom left) of the majorant minimization and exact flux x-component $\frac{\partial u}{\partial x}$ (top right) and exact multiplier $\lambda$ (bottom right).}\label{benchmark_square_discrete_majorant}
\end{figure}

\begin{figure}
\center
\begin{minipage}{0.45\textwidth}
\includegraphics[width=\textwidth]{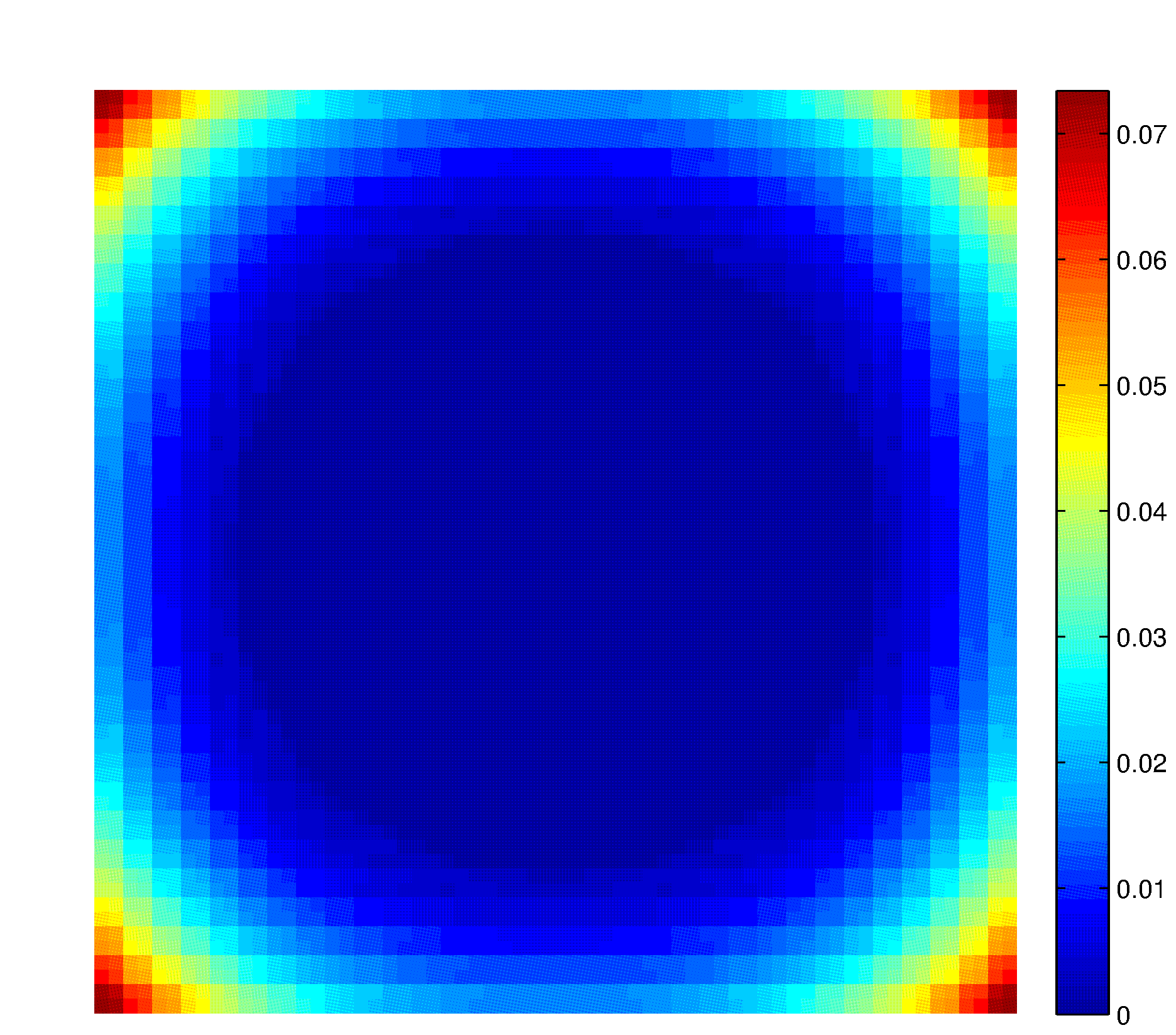}
\end{minipage}
\begin{minipage}{0.45\textwidth}
\includegraphics[width=\textwidth]{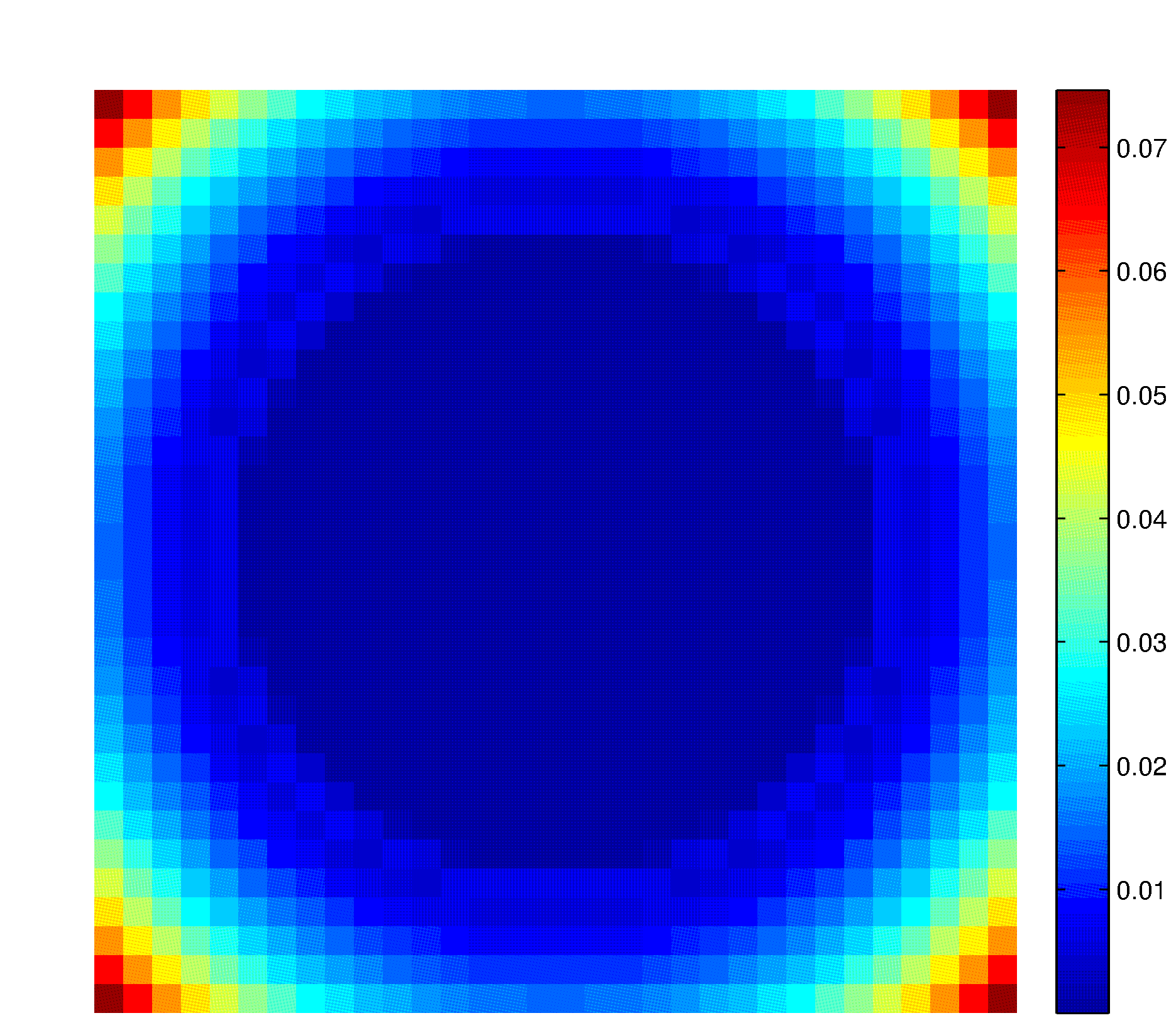}
\end{minipage}
\caption{Local distribution of the exact error $\frac{1}{2}\|v-u\|^{2}_{E}$ 
(left) and of the functional majorant (right).}\label{fig:error_majorant}
\center
\begin{minipage}{0.3\textwidth}
\includegraphics[width=\textwidth]{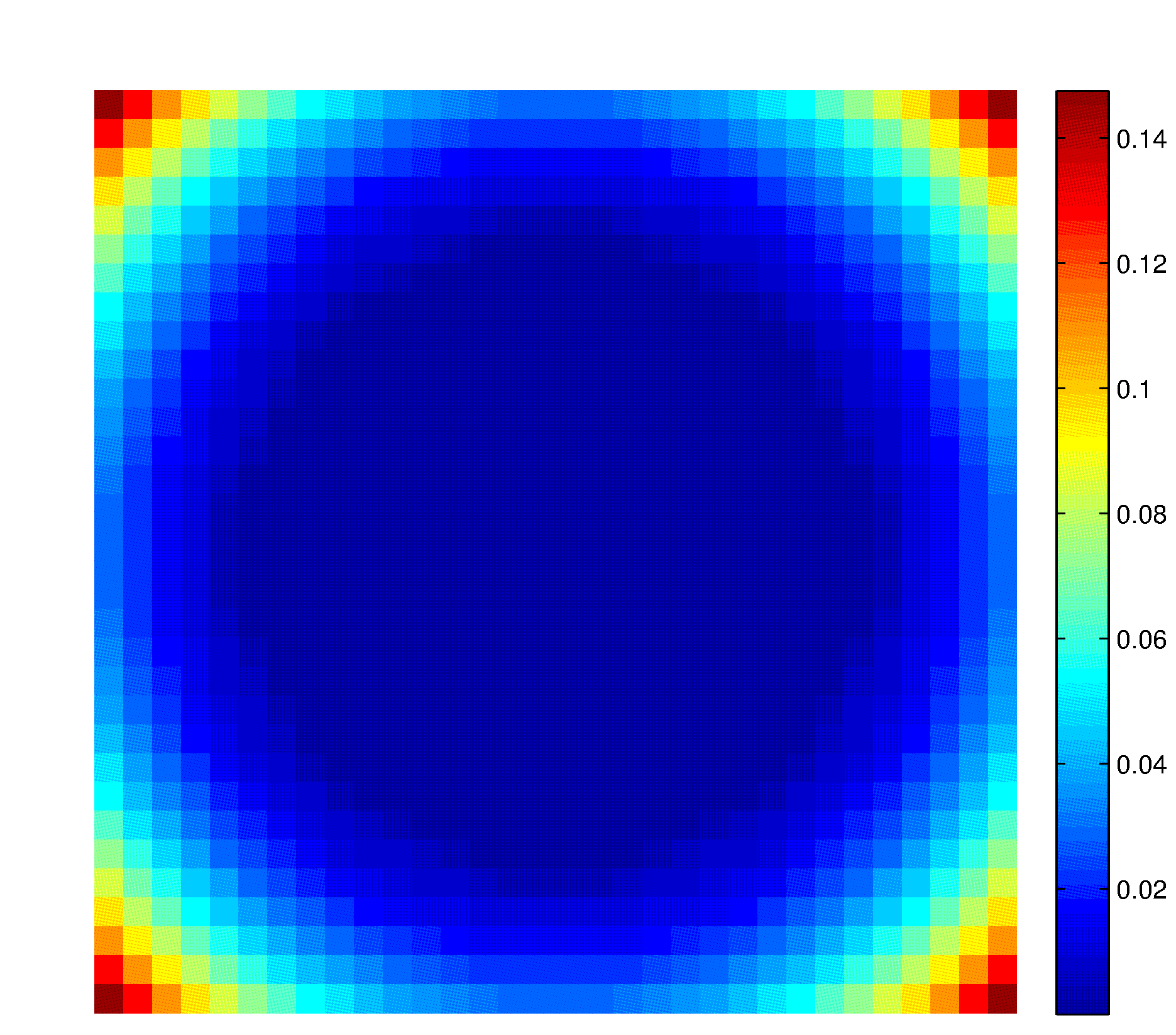}
\end{minipage}
\begin{minipage}{0.32\textwidth}
\includegraphics[width=\textwidth]{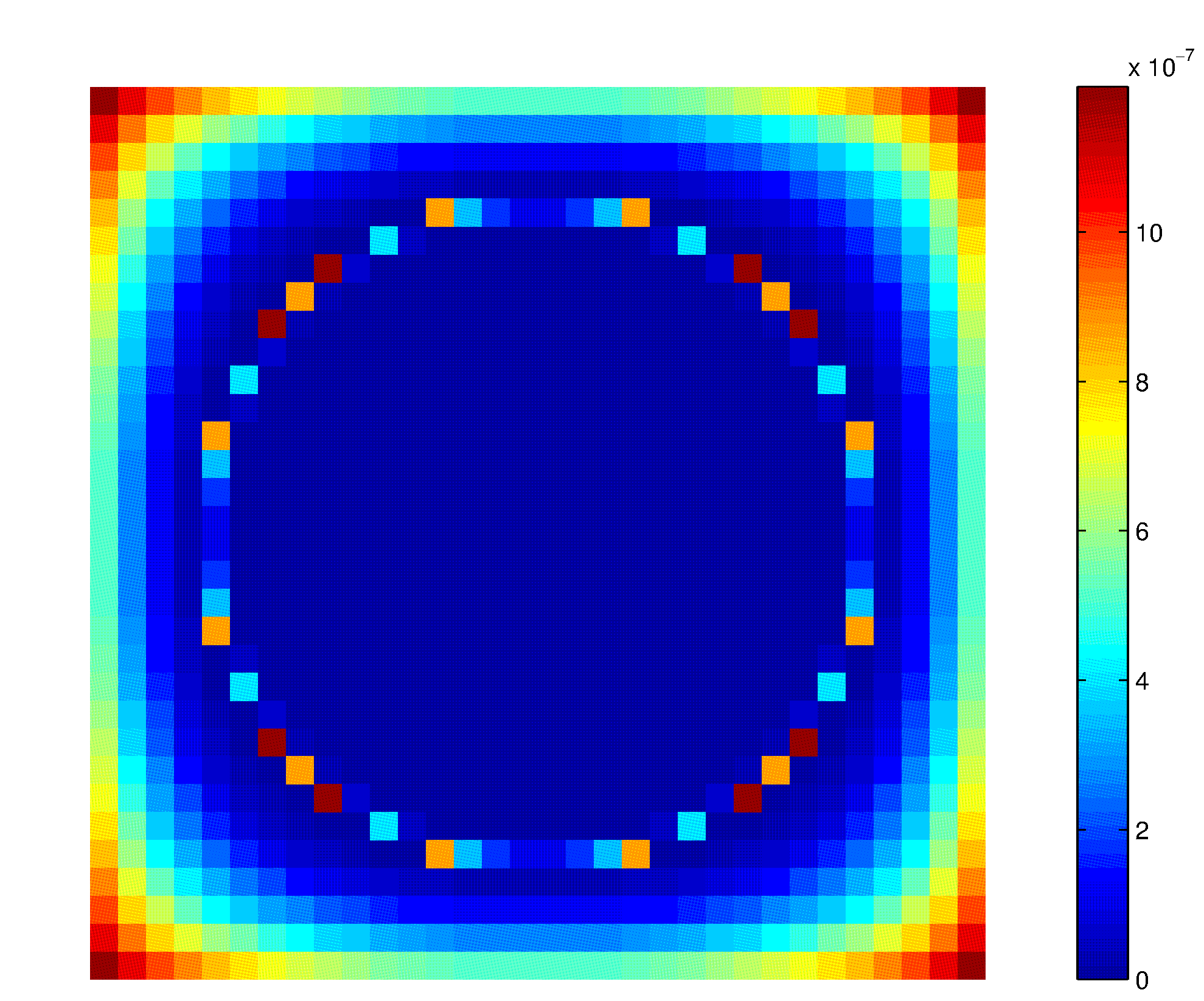}
\end{minipage}
\begin{minipage}{0.3\textwidth}
\includegraphics[width=\textwidth]{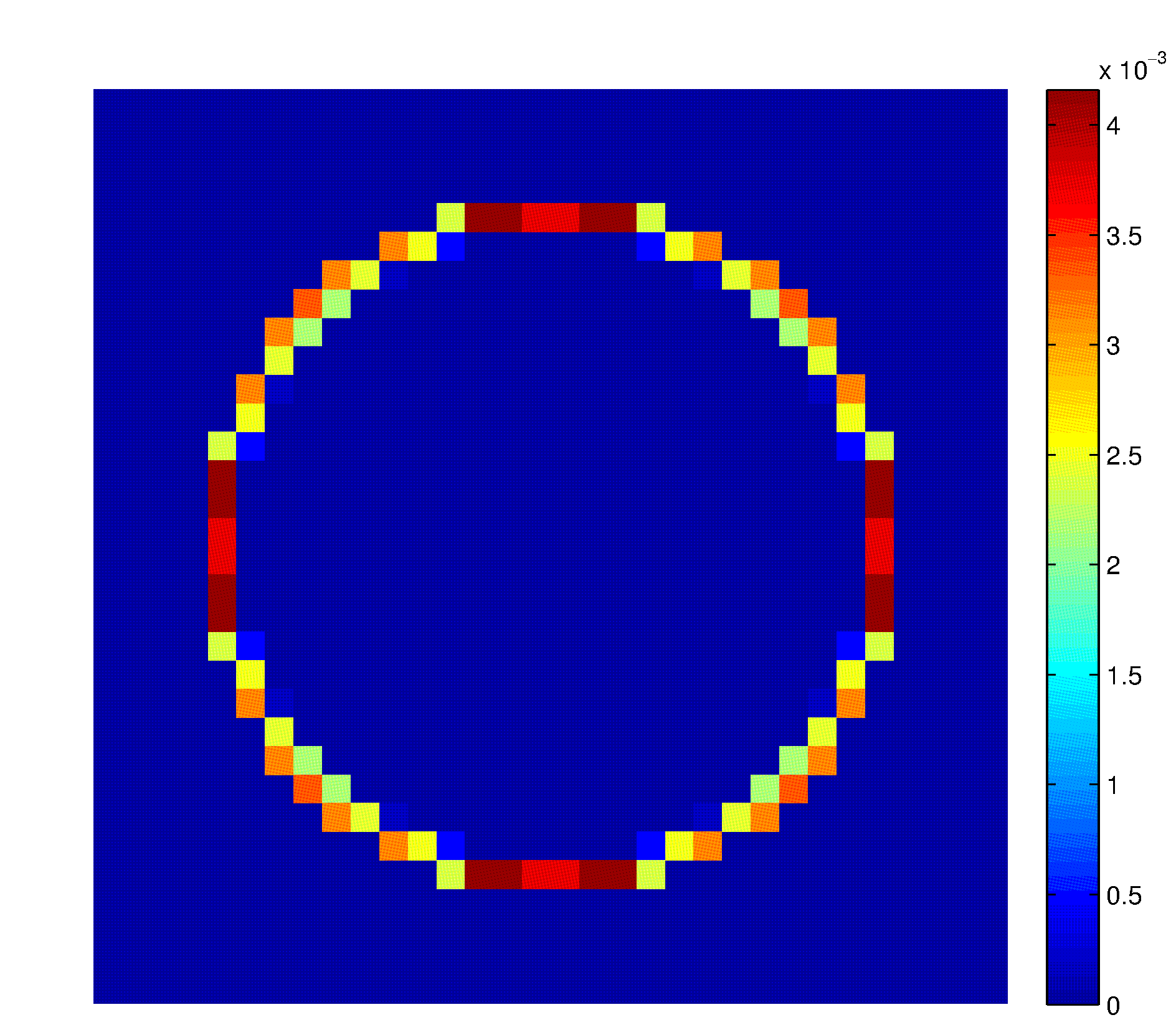}
\end{minipage}
\caption{Local distributions of the functional majorants parts: $\|{\nabla v-\tau^{*}}\|^{2}_\Omega $ (left), $\|{\rm{div}} \,\tau^{*}+f+\mu\|^{2}_\Omega$ (middle), $\intO \mu(v-\phi)\textrm{d}x$ (right). Note that the amplitudes $10^{-7}$ and $10^{-3}$ of the middle and the right pictures are small in comparison with the left picture.
}\label{fig:majorant_parts}
\end{figure}

\begin{figure}
\center
\includegraphics[width=0.8\textwidth]{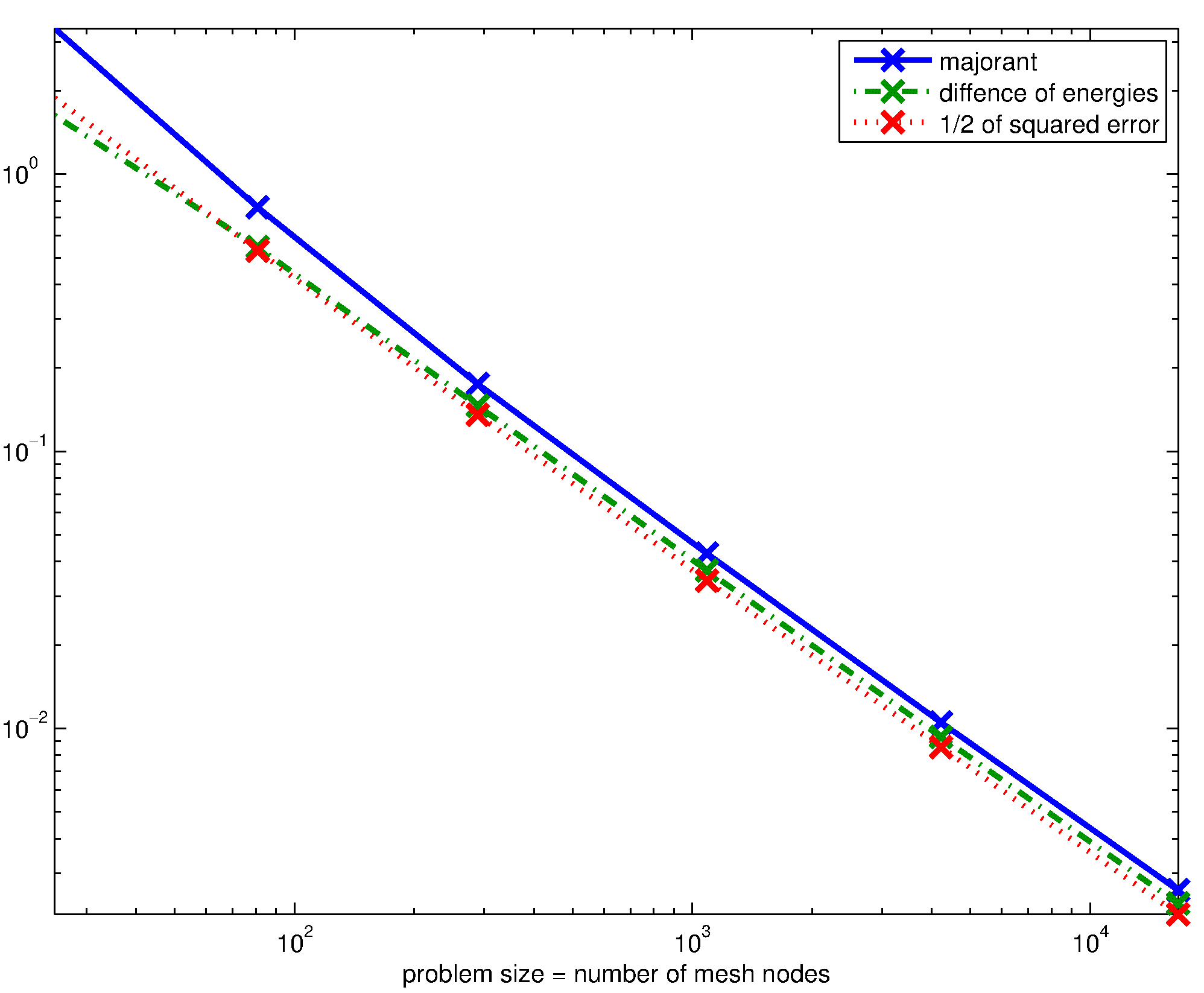}
\caption{Square benchmark with a constant obstacle: convergence.}\label{benchmark_square_convergence}
\end{figure}

In Benchmark I, we consider the case of the contact radius 
$$R=0.7$$
only. The finest considered meshes $\mathcal{T}_{h=1/64}$ and $\mathcal{T}_{h=1/256}$ to provide a discrete solution $v$ and evaluate its error $\|v_h-u\|_{E}$ according to \eqref{eq:error_estimate3} 
are characterized by geometrical properties 
$$ \mathcal{T}_{h=1/64}:  \mbox{ 16384 elements, 16641 nodes and 33024 edges } $$
$$ \mathcal{T}_{h=1/256}: \mbox{ 262144 elements, 263169 nodes and 525312 edges}. $$
Figure \ref{benchmark_square_discrete_solution} displays a discrete solution $v$ computed on a rectangular mesh $\mathcal{T}_{h=1/16}$ and Figure \ref{benchmark_square_discrete_majorant} a discrete flux component $\tau^*_x$ in x-direction and a discrete multiplier $\mu$ computed from the majorant minimization algorithm (Algorithm \ref{Alg:Majorant}) on the same rectangular mesh. 
$\mathcal{T}_{h=1/16}$ was chosen not as the finest rectangular mesh possible in order to stress out the shape of applied finite elements, in particular of discontinuous Raviart-Thomas elements for fluxes. A discrete flux component $\tau^*_y$ in y-direction is not shown due to symmetry with $\tau^*_x$. 
Local distributions of the exact error and of the functional majorant are compared in Figure \ref{fig:error_majorant} and three majorant parts are further depicted in Figure \ref{fig:majorant_parts}. 
We notice that the majorant part $\|{\nabla v-\tau^{*}}\|^{2}_\Omega $ has significantly higher amplitude than the equilibrium part $\|{\rm{div}} \,\tau^{*}+f+\mu\|^{2}_\Omega$ and the nonlinear part $\intO \mu(v-\phi)\textrm{d}x$. The nonlinear part is further nonzero distributed around a thin layer boundary of the contact domain corresponding to the discrete solution $v$. 
Converge behaviour for all considered rectangulations is compared in Figure \ref{benchmark_square_convergence}. 
We can see that the energy estimate \eqref{energyEstimate} and the majorant estimate \eqref{majorantEstimate} are very sharp with valid inequalities signs. 

\begin{remark}
A simplified least-square (LS) variant of the function majorant type estimate was tested on the same benchmark including a mesh adaptivity in \cite{CaMe}. The main difference is that nonlinear term 
$\intO \mu(v-\phi)\textrm{d}x$ is considered here in order to quarantee the majorant upper bound in \eqref{majorantEstimate} and enhance an accurate error control. 
\end{remark}

\begin{figure}
\center
\includegraphics[width=0.4\textwidth]{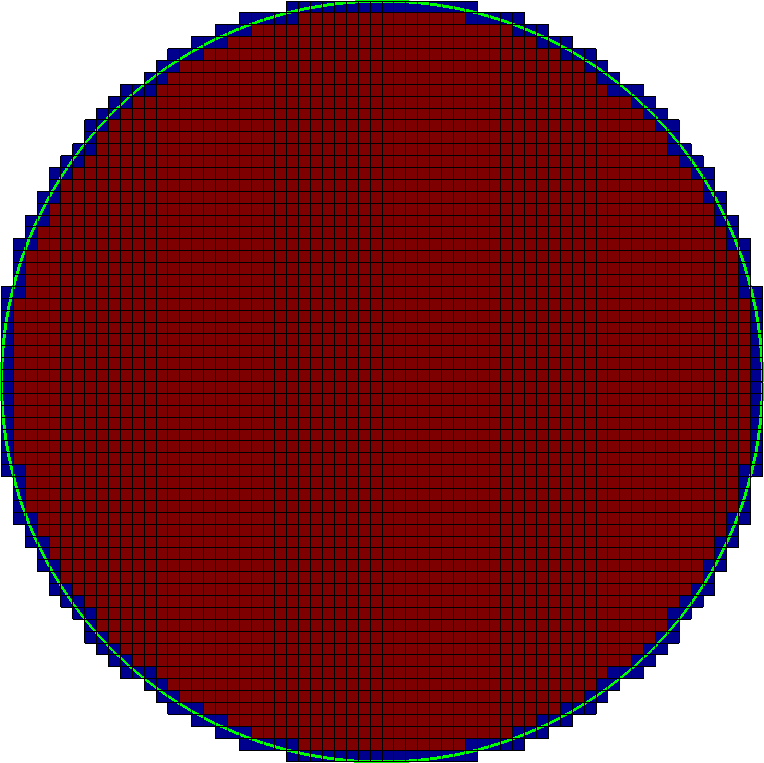}
\caption{Rectangulation of the ring domain. The green line indicates the exact ring boundary, red rectangles are completely inscribed in the ring boundary and blue circumscribed rectangles contain at least one node lying outside the ring boundary. Red rectangles define an inscribed rectangulation $\mathcal{T}_\vee$ and red and blue rectangles together define a circumscribed rectangulation $\mathcal{T}^\wedge$.} 
\label{fig:ring}
\end{figure}

\begin{figure}
\center
\begin{minipage}{0.49\textwidth}
\center
\includegraphics[width=0.8\textwidth]{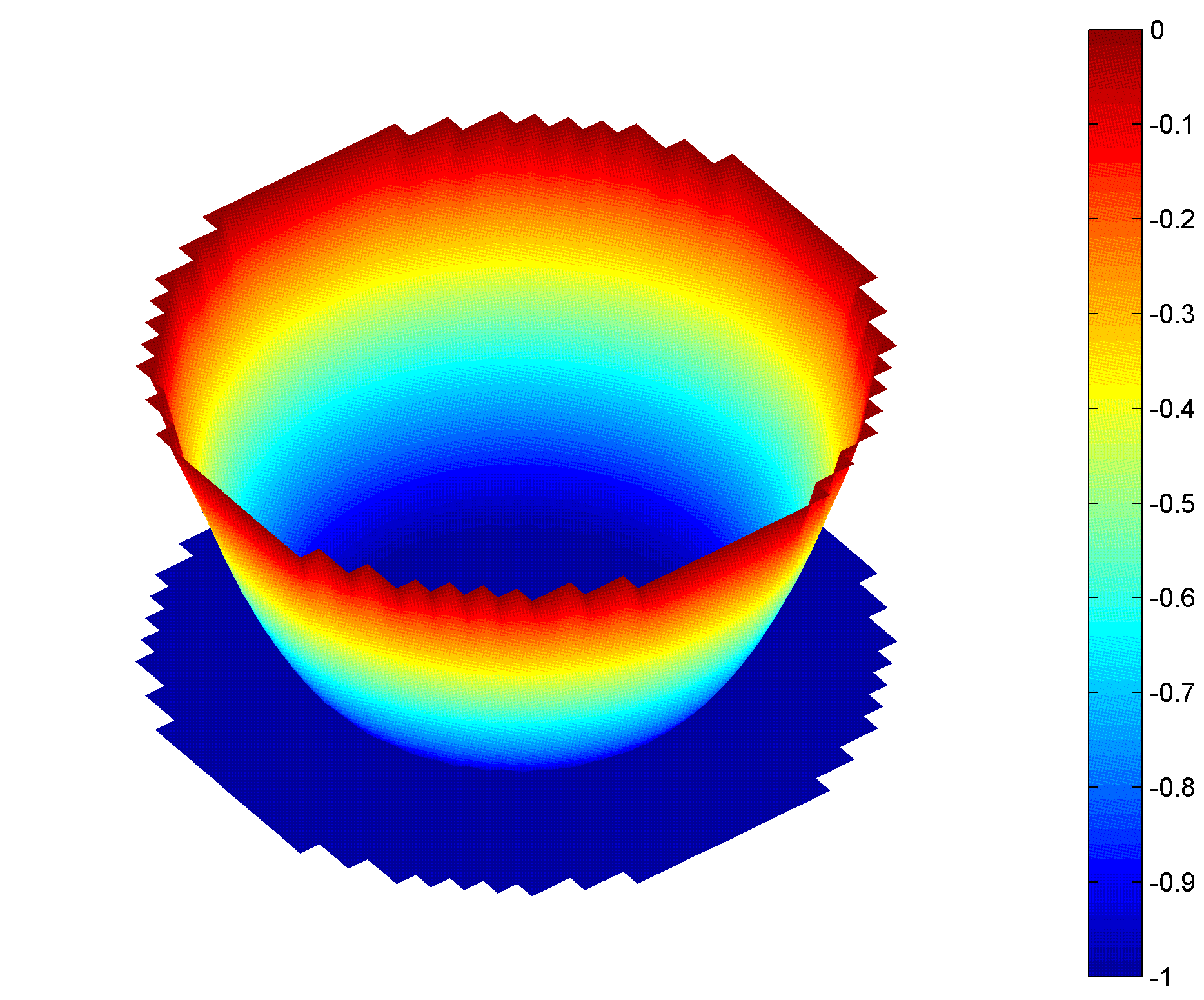}
\end{minipage}
\begin{minipage}{0.49\textwidth}
\center
\includegraphics[width=0.8\textwidth]{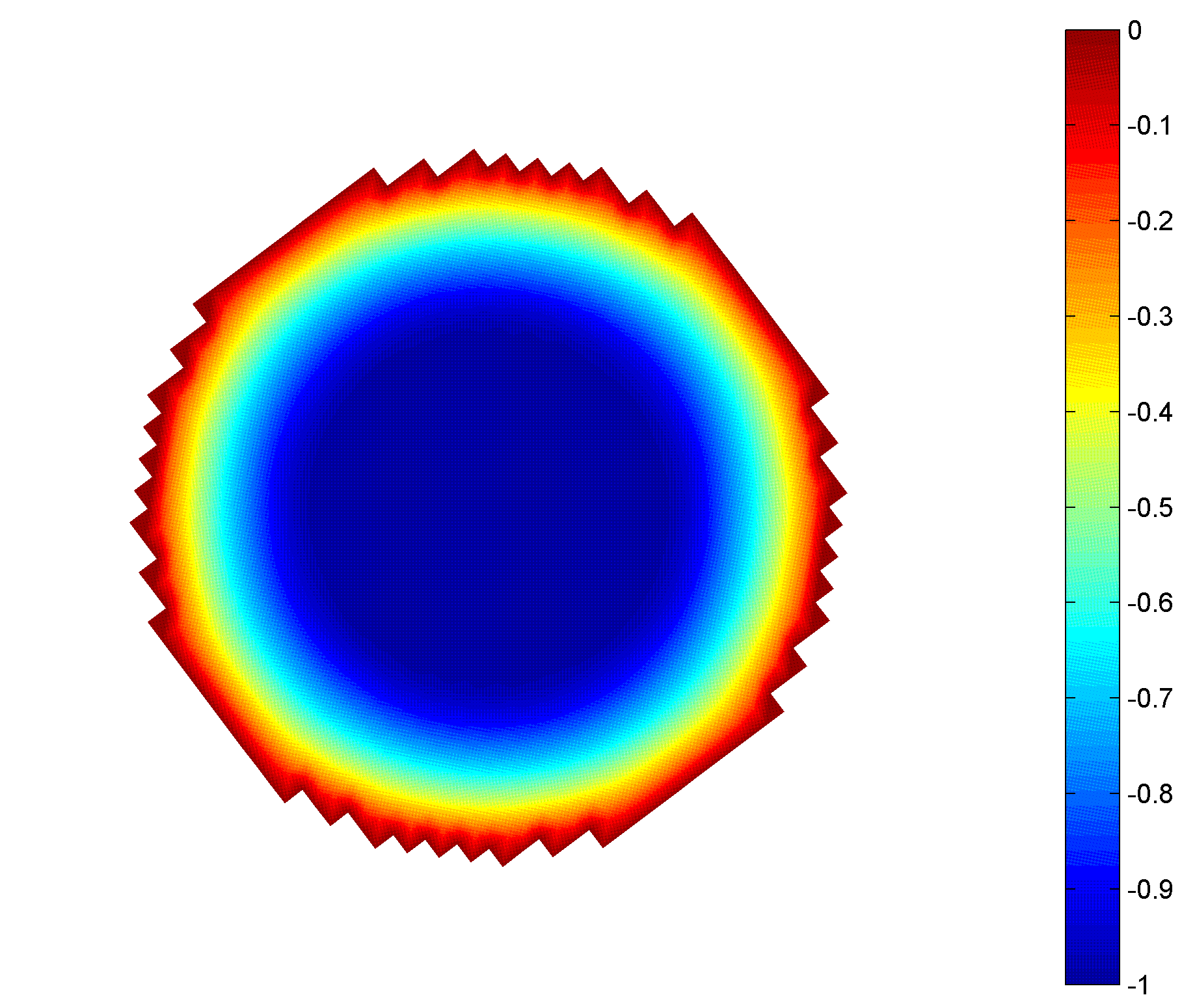}
\end{minipage}
\caption{Discrete solution $v$ of the obstacle problem and the lower obstacle $\phi$ (left) and its rotated view (right). The dark blue color indicates the contact domain.}
\label{benchmark_ringConstantObstacle_discrete_solution}
\begin{minipage}{0.49\textwidth}
\center
\includegraphics[width=0.8\textwidth]{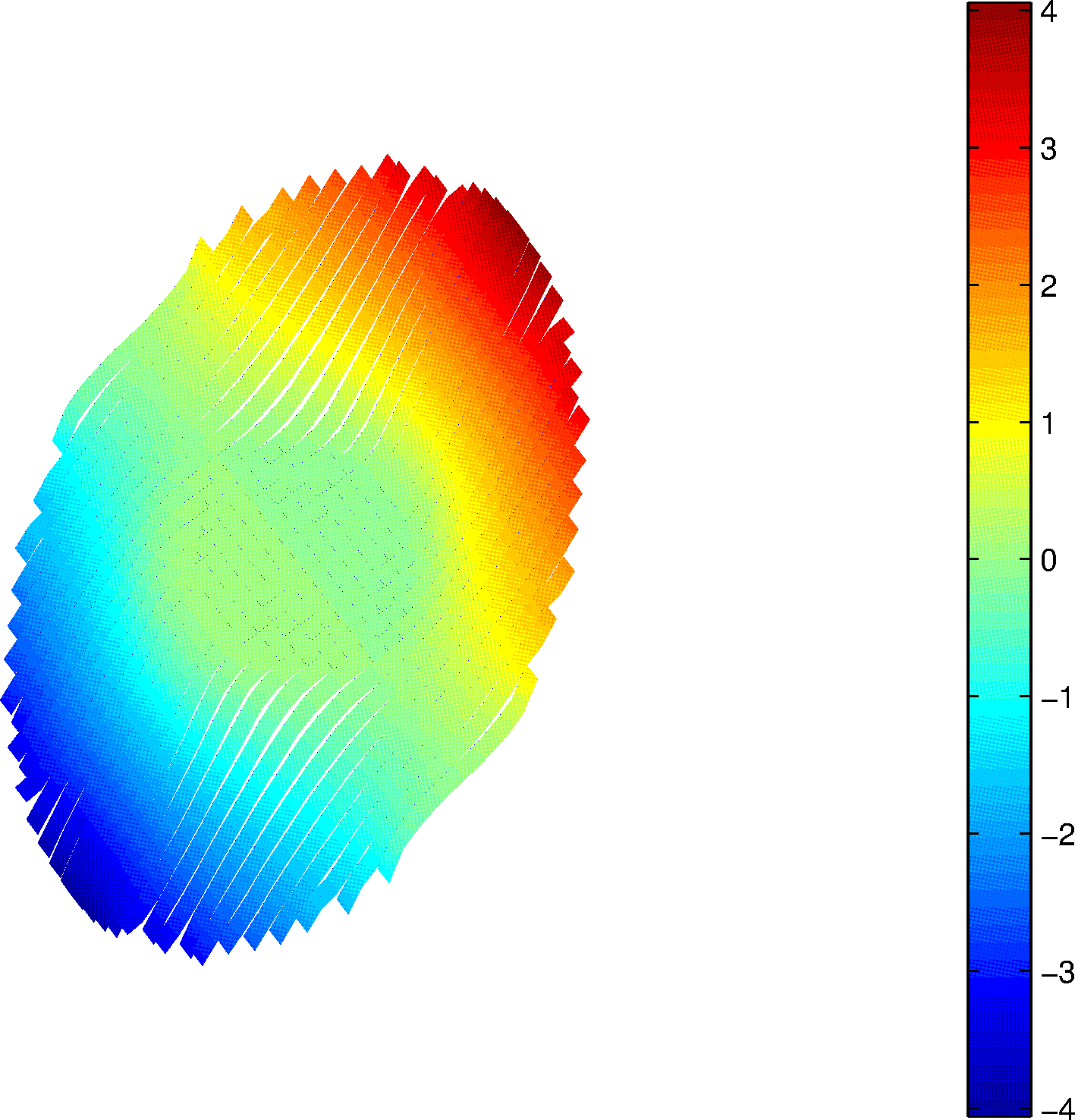}
\vspace{0.2cm} \\
\includegraphics[width=0.8\textwidth]{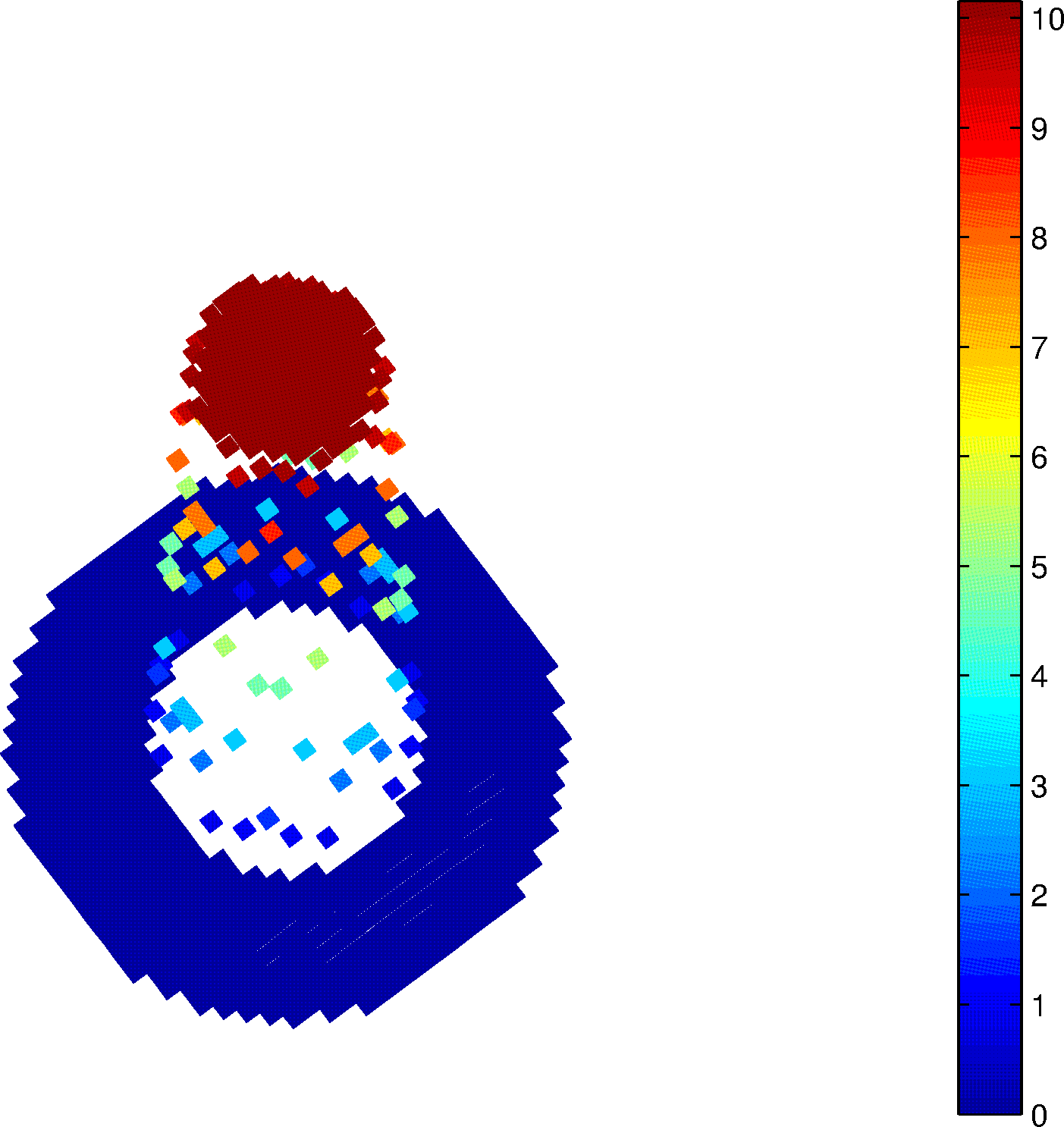}
\end{minipage}
\begin{minipage}{0.49\textwidth}
\center
\includegraphics[width=0.8\textwidth]{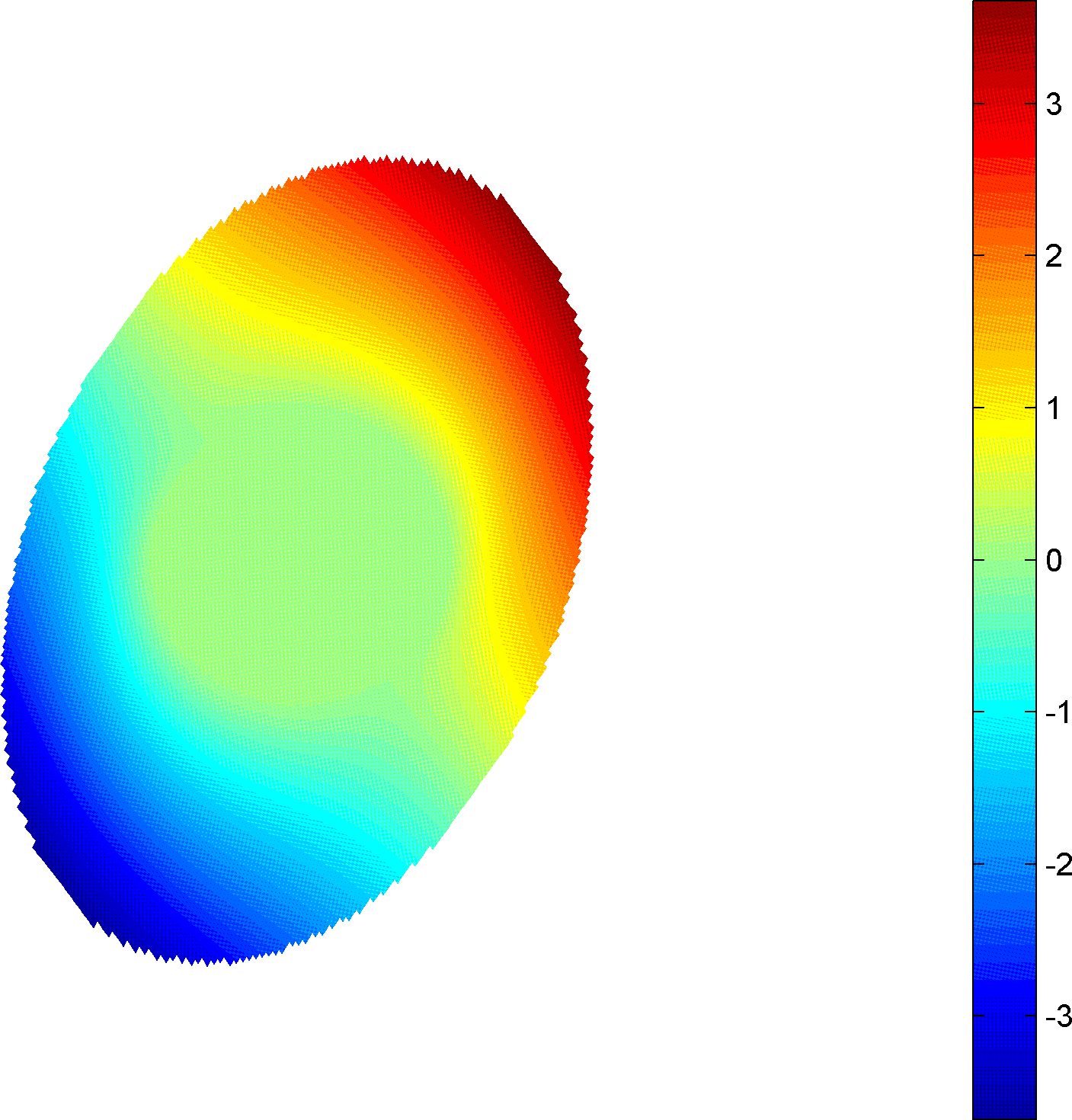} 
\vspace{0.2cm} \\
\includegraphics[width=0.8\textwidth]{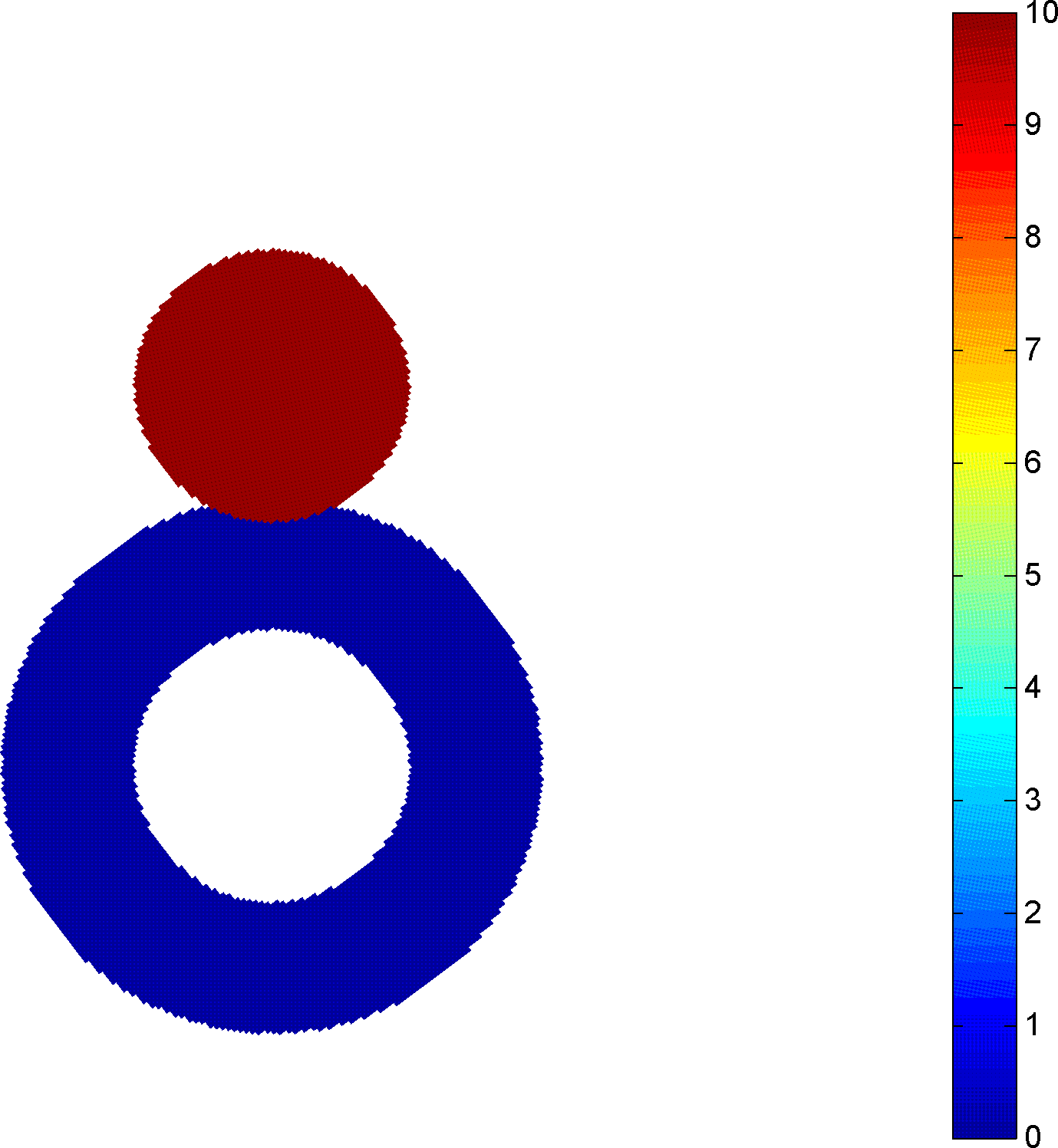}
\end{minipage}
\caption{Discrete flux x-component $\tau_x^*$ (top left) and discrete multiplier $\mu$ (bottom left) of the majorant minimization and exact flux x-component $\frac{\partial u}{\partial x}$ (top right) and exact multiplier $\lambda$ (bottom right).}\label{benchmark_ringConstantObstacle_discrete_majorant}
\end{figure}

\begin{figure}
\center
\begin{minipage}{0.49\textwidth}
\center
\includegraphics[width=0.8\textwidth]{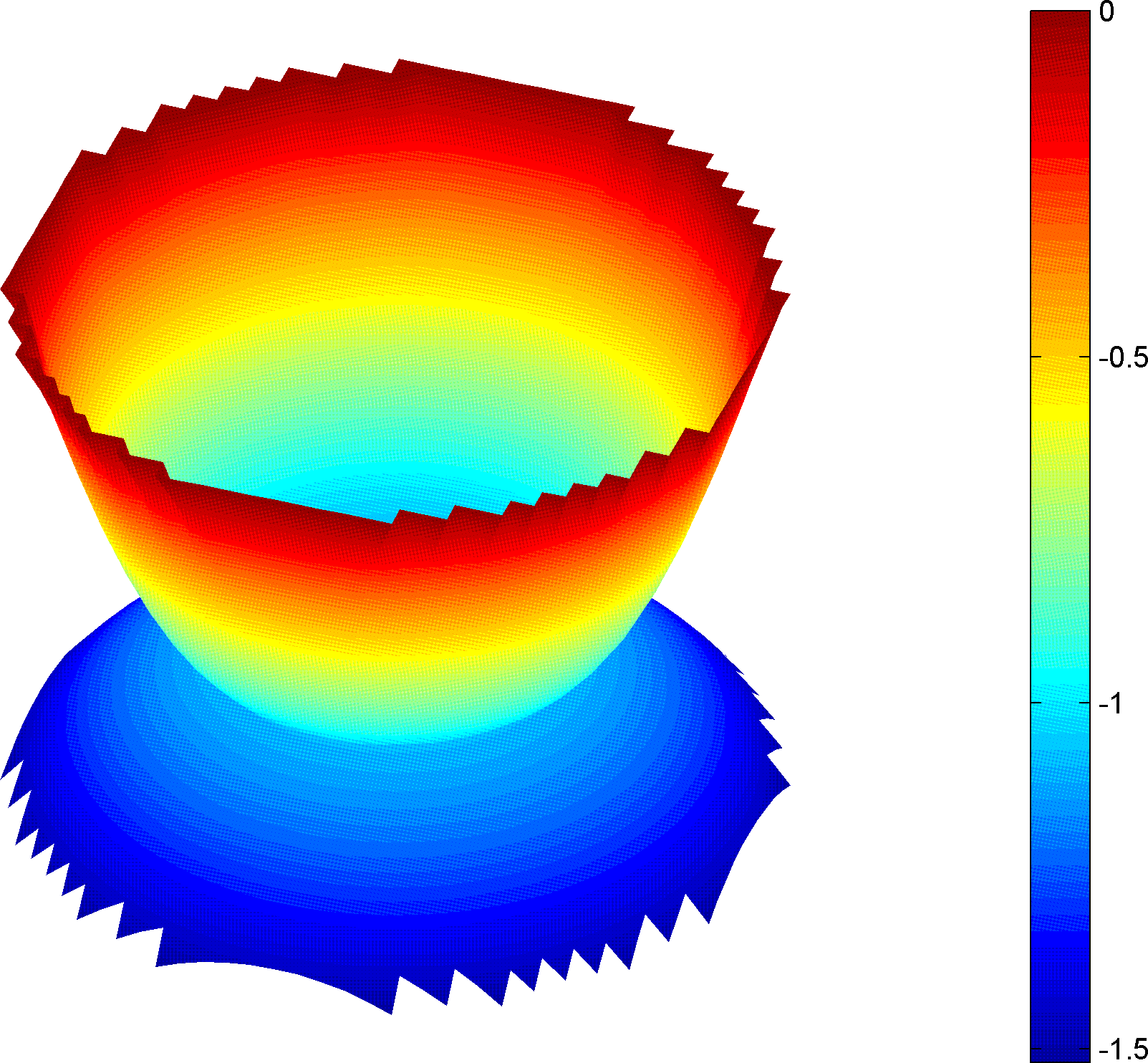}
\end{minipage}
\begin{minipage}{0.49\textwidth}
\center
\includegraphics[width=0.8\textwidth]{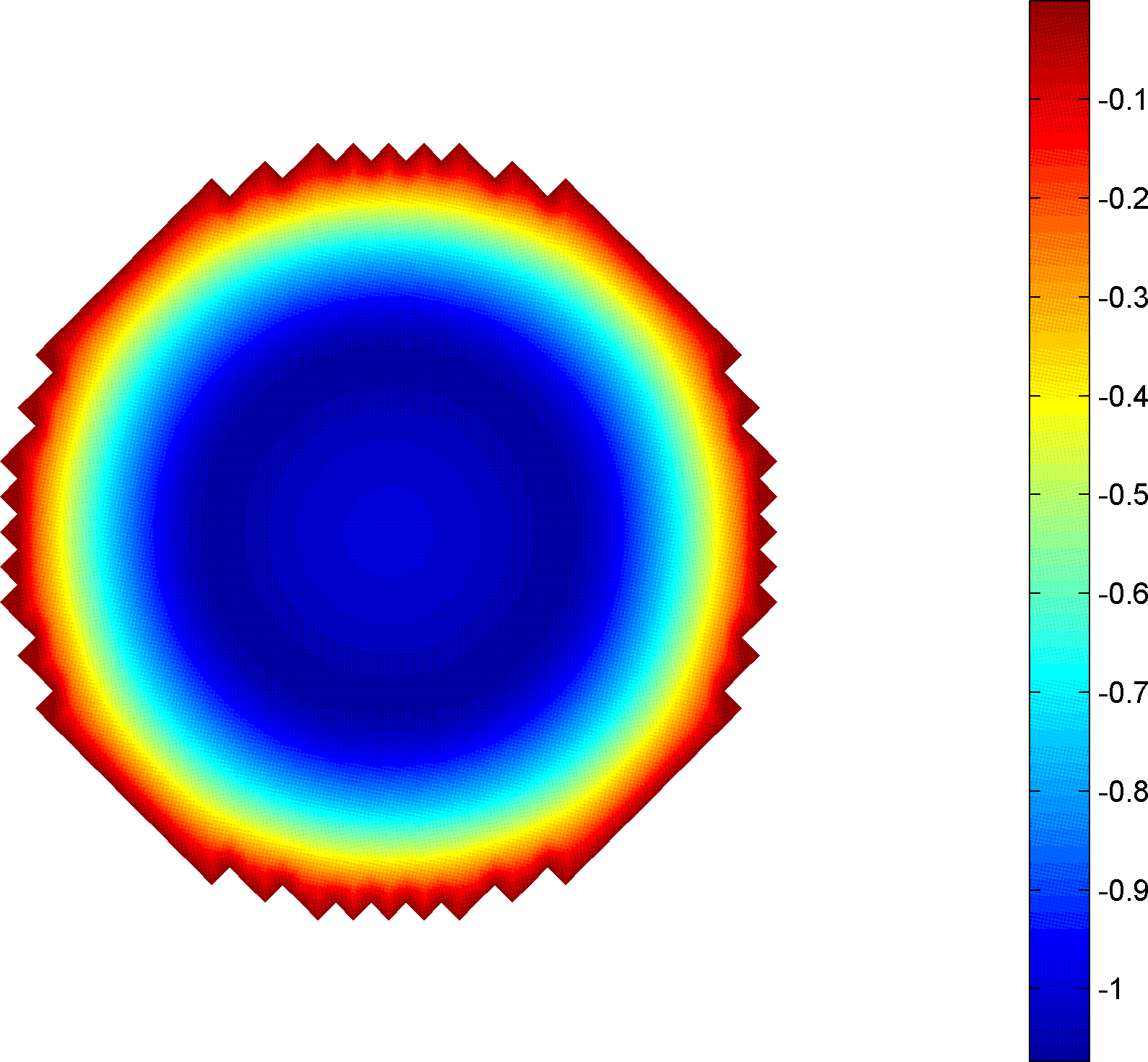}
\end{minipage}
\caption{Discrete solution $v$ of the obstacle problem and the lower obstacle $\phi$ (left) and its rotated view (right). The dark blue color indicates the contact domain.}
\label{benchmark_ringSphericalObstacle_discrete_solution}
\begin{minipage}{0.49\textwidth}
\center
\includegraphics[width=0.8\textwidth]{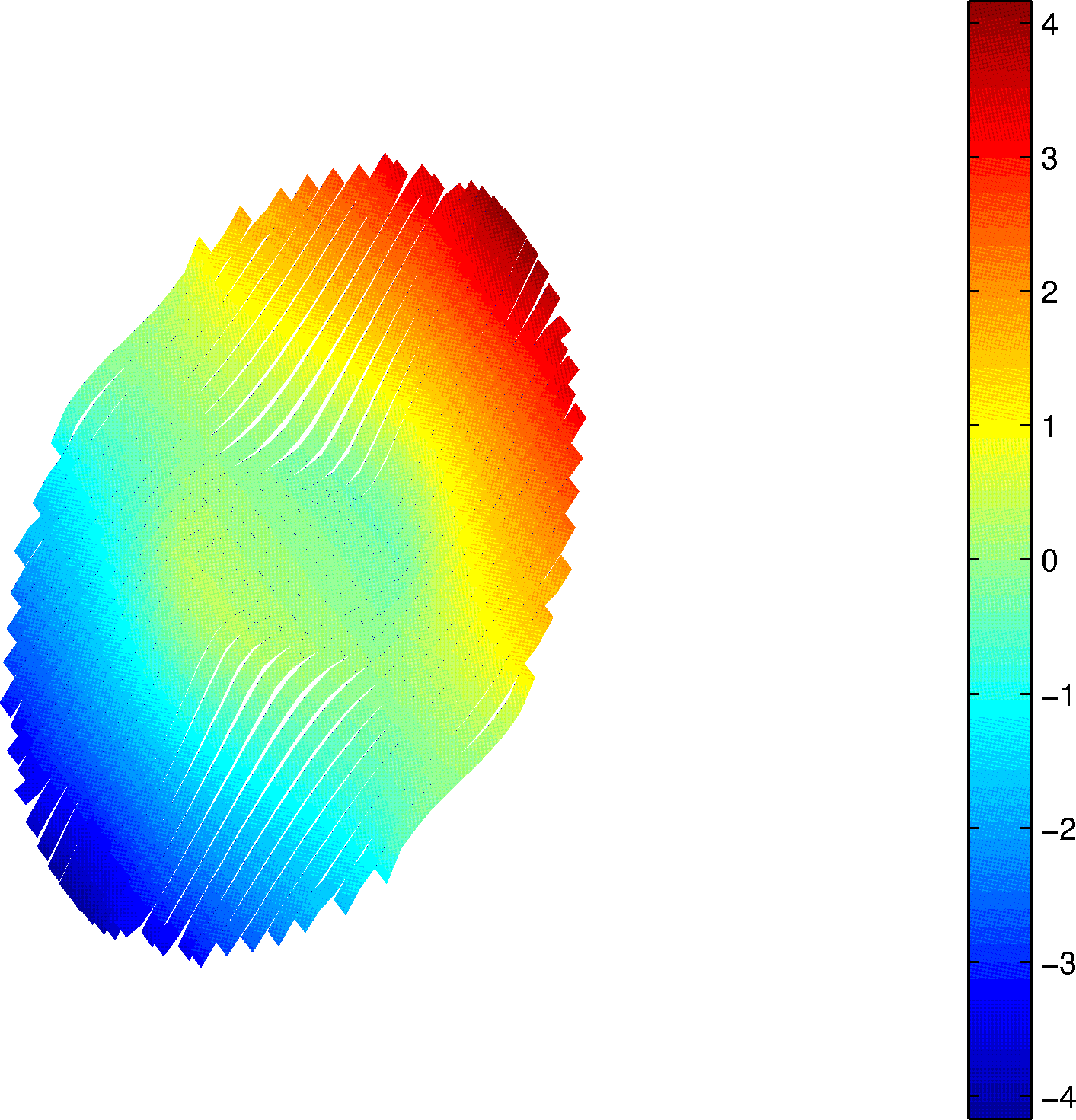}
\vspace{0.2cm} \\
\includegraphics[width=0.8\textwidth]{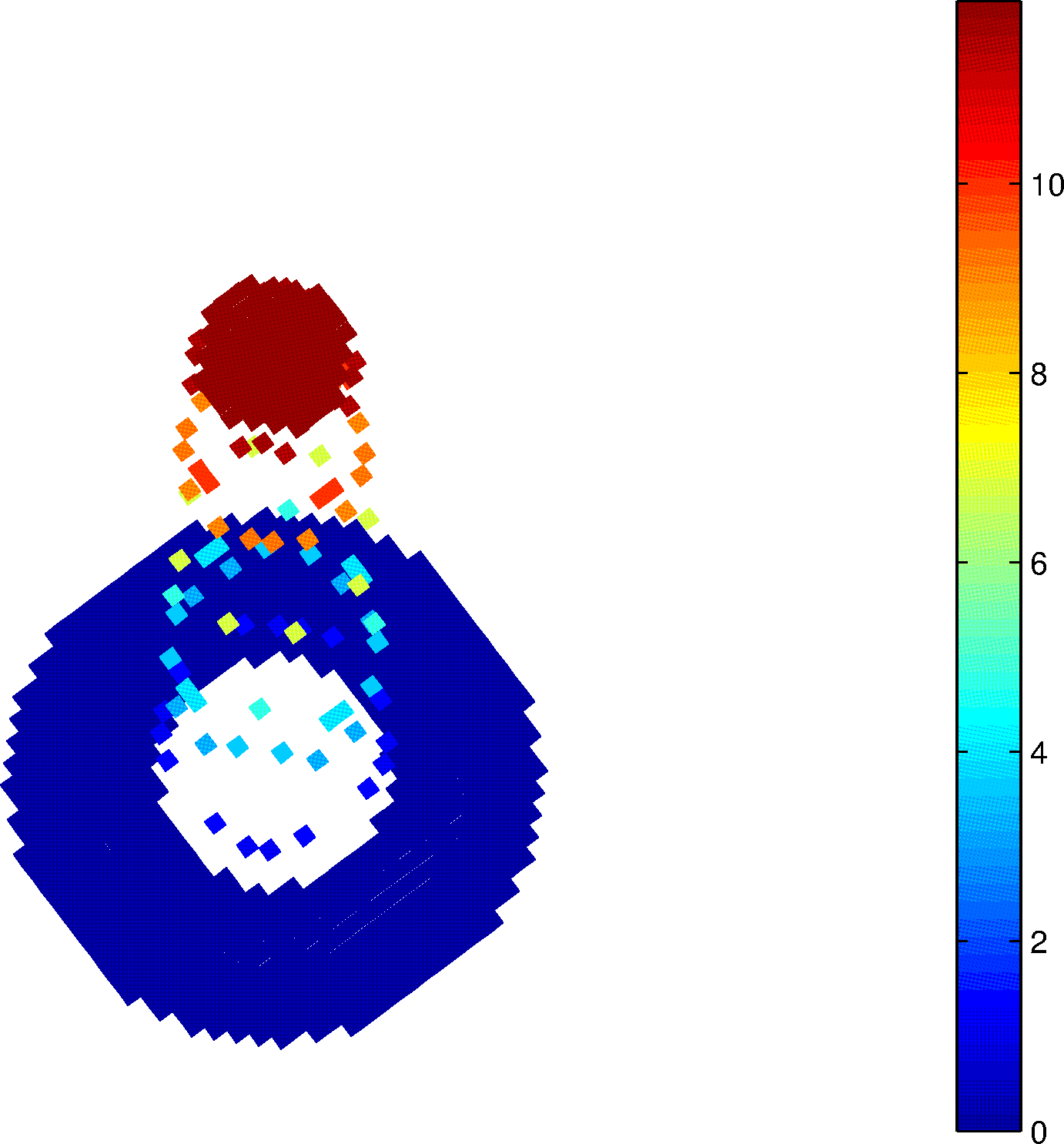}
\end{minipage}
\begin{minipage}{0.48\textwidth}
\center
\includegraphics[width=0.8\textwidth]{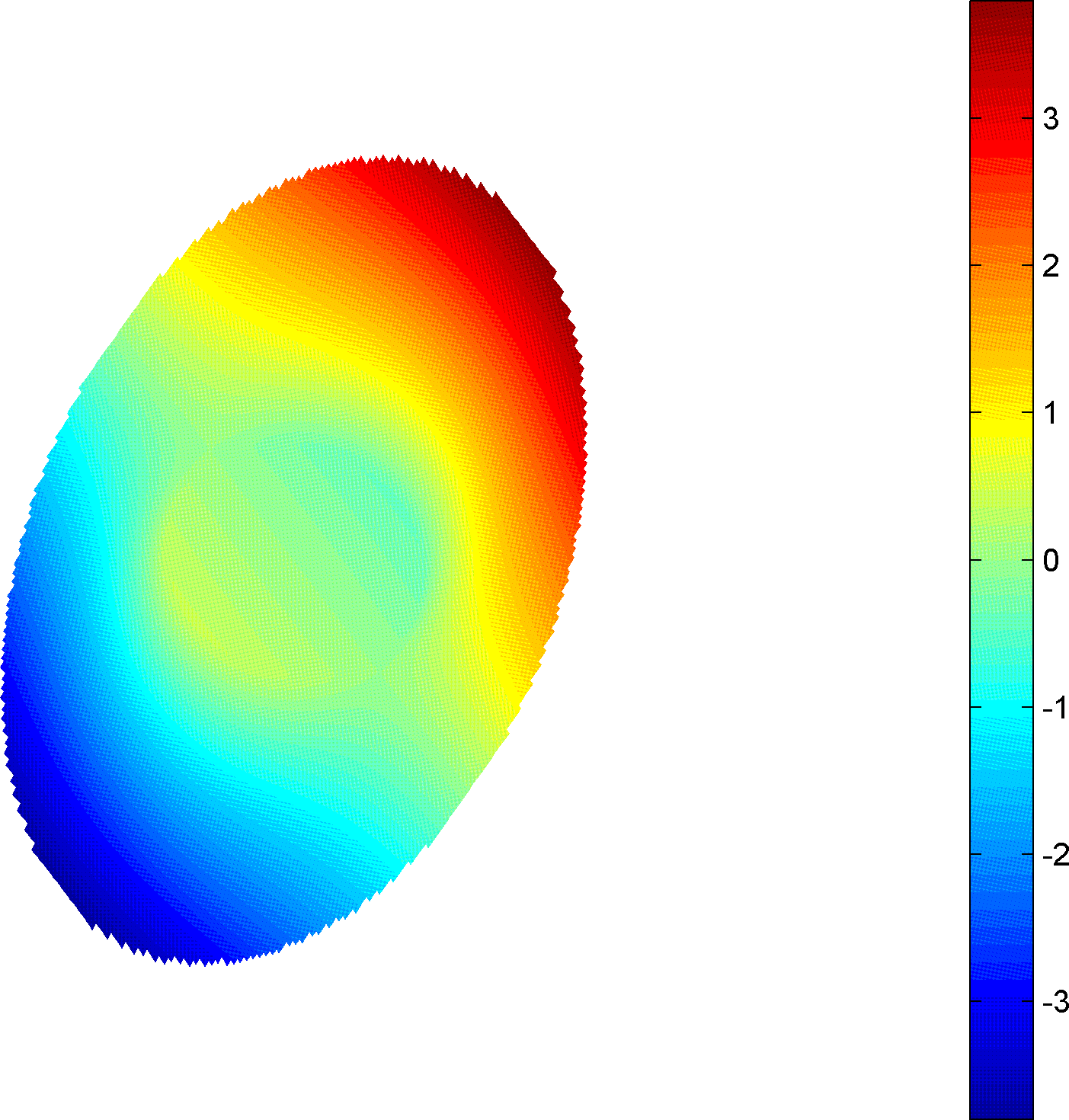}
\vspace{0.2cm} \\
\includegraphics[width=0.8\textwidth]{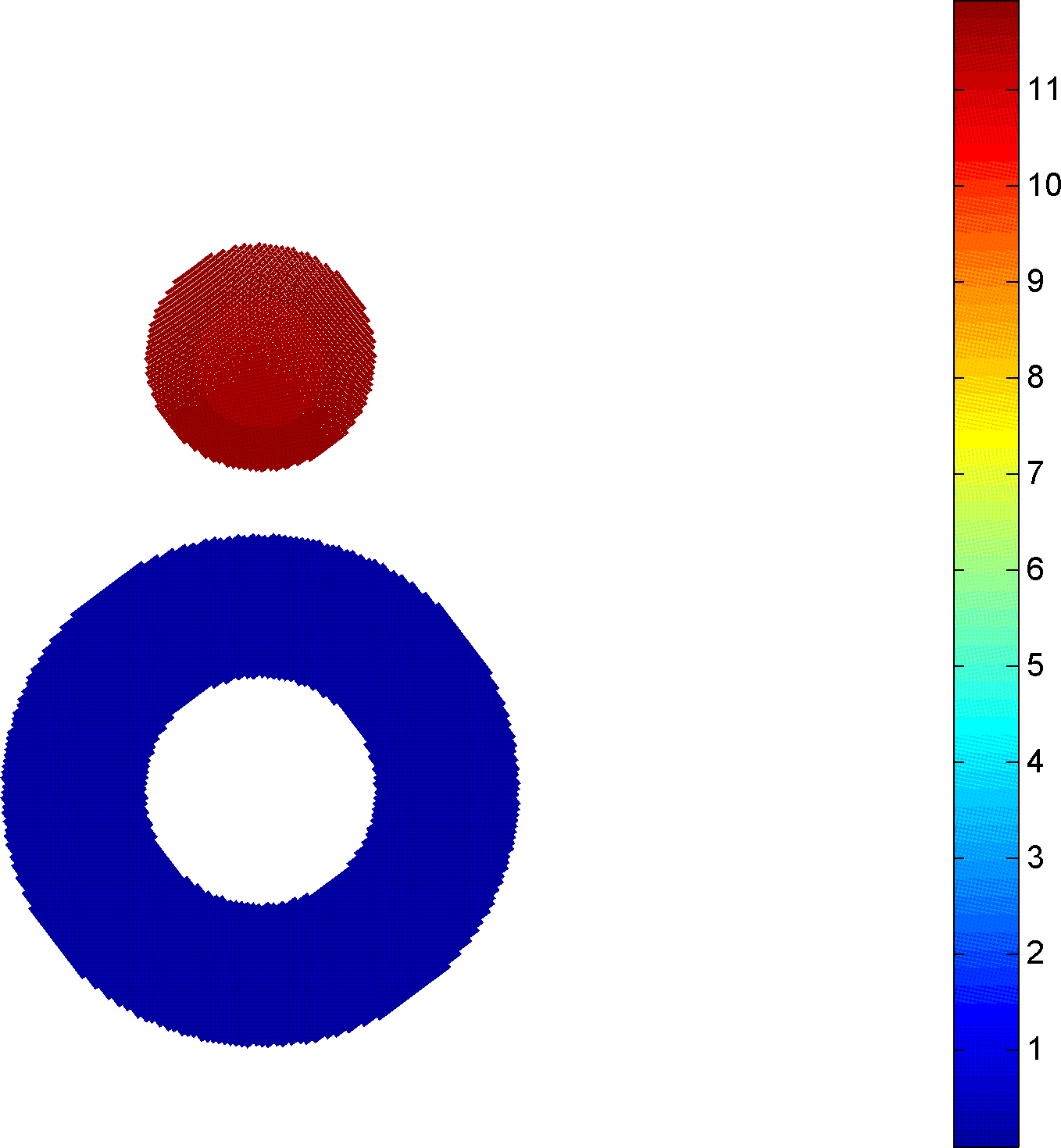}
\end{minipage}
\caption{Discrete flux x-component $\tau_x^*$ (top left) and discrete multiplier $\mu$ (bottom left) of the majorant minimization and exact flux x-component $\frac{\partial u}{\partial x}$ (top right) and exact multiplier $\lambda$ (bottom right).}\label{benchmark_ringSphericalObstacle_discrete_majorant}
\end{figure}

\begin{figure}
\center
\includegraphics[width=0.7\textwidth]{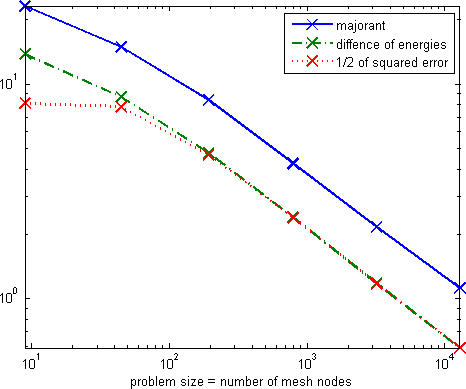}
\caption{Ring benchmark with a constant obstacle: convergence.}\label{benchmark_ringConstantObstacle_convergence}
\end{figure}

\begin{figure}
\center
\includegraphics[width=0.7\textwidth]{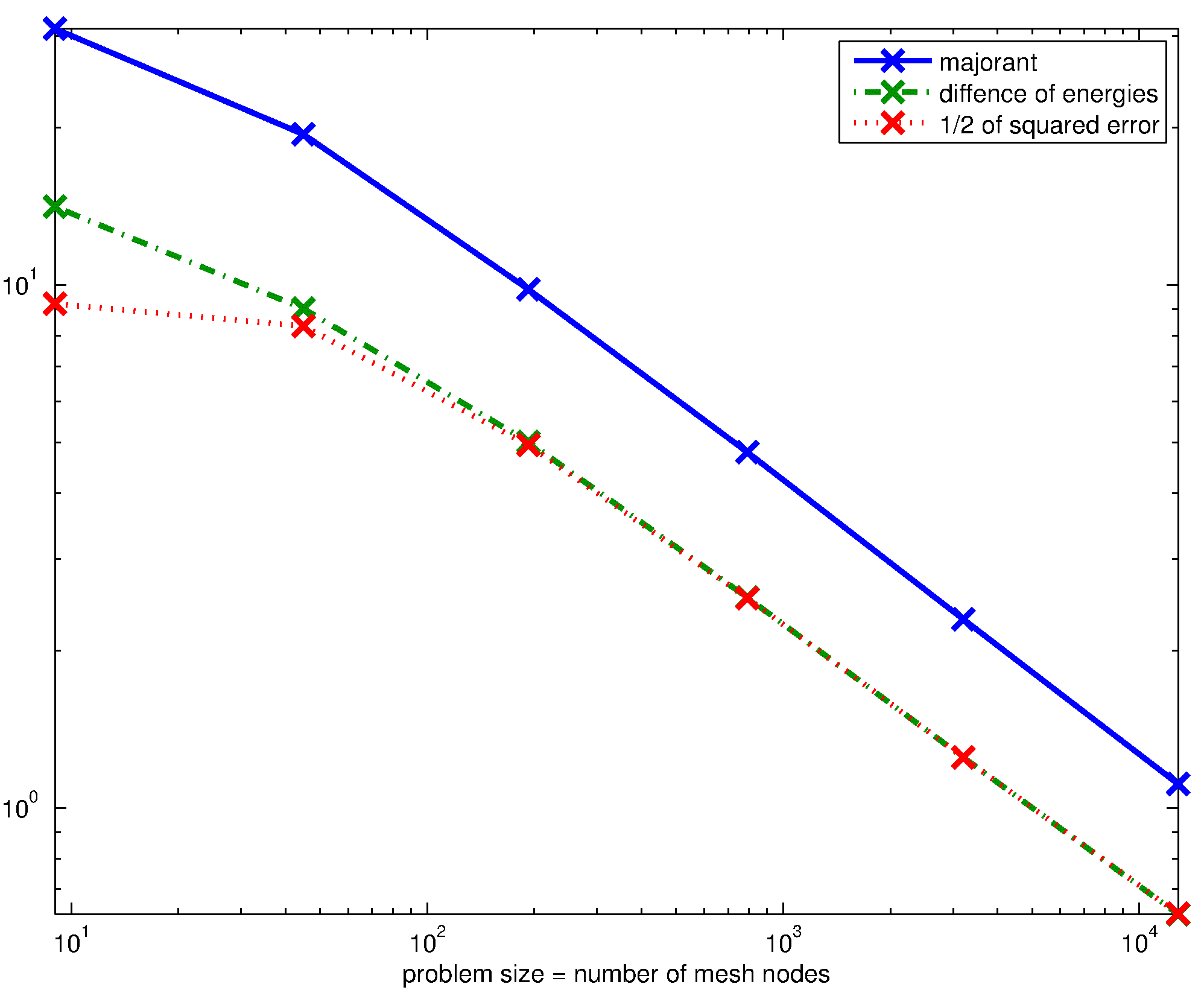}
\caption{Ring benchmark with a spherical obstacle: convergence. } \label{benchmark_ringSphericalObstacle_convergence}
\end{figure}

In Benchmark II and Benchamark III, we consider the cases
$$ f=-10, \phi=-1 \mbox{ (for Benchmark II)} $$
$$ f=-10, \phi_{\mathrm{max}}=-1, \rho=1.2 \mbox{ (for Benchmark III)}$$
only. Numerical solutions of \eqref{condR_c} and \eqref{condR_s} show that contact radius $R$ is approximately
$R \approx  0.5024744$ for Benchmark II and $R \approx  0.4389205$ for Benchmark III.
Since our implementation runs on rectangular elements allowing a polygonal boundary only, 
we consider an inscribed rectangulation  $\mathcal{T}_\vee$ and a circumscribed rectangulation $\mathcal{T}^\wedge$ as approximations of the ring boundary, see Figure \ref{fig:ring} for details. 
A discrete solution $v_\vee$ is solved on 
$\mathcal{T}_\vee$ satisfying zero Dirichlet boundary conditions on 
its 
boundary $\partial \mathcal{T}_\vee$. Then, $v_\vee$ is extended by zero values on $\mathcal{T}^\wedge \setminus \mathcal{T}_\vee$ (displayed by the blue color rectangles in Figure \ref{fig:ring}) to define a discrete solution $v^\wedge$ on the circumscribed rectangulations $\mathcal{T}^\wedge$.
It can be easily checked that 
$$J(v_\vee) := \frac{1}{2}\int_{\mathcal{T}_\vee}\nabla v_\vee \cdot \nabla v_\vee \,\textrm{d}x-\int_{\mathcal{T}_\vee} f_\vee v_\vee \,\textrm{d}x 
             = \frac{1}{2}\int_{\mathcal{T}^\wedge}\nabla v^\wedge \cdot \nabla v^\wedge \,\textrm{d}x-\int_{\mathcal{T}^\wedge} f^\wedge v^\wedge \,\textrm{d}x:= J(v^\wedge), $$
where $f_\vee$ represents a restriction of $f$ to $\mathcal{T}_\vee$ and $f^\wedge$ represents an extension of $f_\vee$ to $\mathcal{T}^\wedge \setminus \mathcal{T}_\vee$ by any value. 
The extension $f^\wedge$ is defined by the same constant function $f$ in our implementation. Finally, the majorant minimization is computed on $\mathcal{T}^\wedge$. 
Thus the modification of the energy estimate \eqref{energyEstimate} and the majorant estimate \eqref{majorantEstimate} is combined in the estimate
\begin{equation}\label{energyAndMajorantEstimate_modified} 
\frac{1}{2}\big|\!\big| v_\vee-u|_{\mathcal{T}_\vee} \big|\!\big|^{2}_{E} \leq 
J(v^\wedge)-J(u)\leq \mathcal{M}(v^\wedge,f^\wedge,\phi^\wedge;\beta,\mu^\wedge,\tau^{*\wedge}),
\end{equation}
which is shown in convergence figures. Figures \ref{benchmark_ringConstantObstacle_discrete_solution} and \ref{benchmark_ringSphericalObstacle_discrete_solution} display discrete solutions $v$ of Benchmark II and Benchmark III computed on rectangulation created for $h=\frac{1}{16}$ and Figures \ref{benchmark_ringConstantObstacle_discrete_majorant} and \ref{benchmark_ringSphericalObstacle_discrete_majorant}  a discrete flux component in x-direction $\tau_x$ and a  discrete multiplier $\mu$ computed from the majorant minimization algorithm (Algorithm \ref{Alg:Majorant}) on the same rectangulation. Converge behaviour for all considered rectangulations is compared in Figures  \ref{benchmark_ringConstantObstacle_convergence} and \ref{benchmark_ringSphericalObstacle_convergence}. 
We can see that the modified energy and majorant estimates \eqref{energyAndMajorantEstimate_modified} are sharp. 

%

%


\section*{Conclusions}
Computations for three discussed benchmarks with known analytical solutions demonstrate that the functional majorant serves as a fully computable tool to estimate the upper bound of the difference of energies $J(v)-J(u)$ which serves further as an upper bound of the error in the energy norm. The majorant minimization algorithm consists of the solution of a linear system of equations for a flux variable and elementwise computation of the Lagrange multiplier. If a good a initial Lagrange multiplier field is available together with the discrete solution of the obstacle problem, the majorant minimization algorithm requires only few iterations to provide a sharp upper bound of the error. 

\section*{Appendix - local FEM matrices}\label{sec:appendix}
We assume a reference rectangle $T_{\rm{ref}}$ with lengths $h_x$ and $h_y$ defined by vertices 
$$ v_1=(0, 0), \qquad v_2=(0, h_x), \qquad v_3=(h_x, h_y), \qquad v_4=(0, h_y) $$ 
and define four local bilinear nodal basic functions 
\begin{eqnarray*}
&&\hat \psi_1(x,y)=1-\frac{y}{h_y}-\frac{x}{h_x}+\frac{xy}{h_{x}h_{y}},  \\
&&\hat \psi_2(x,y)=\frac{x}{h_x}-\frac{xy}{h_{x}h_{y}}, \quad \\
&&\hat \psi_3(x,y)=\frac{xy}{h_{x}h_{y}}, \\
&&\hat \psi_4(x,y)=\frac{y}{h_y}-\frac{xy}{h_{x}h_{y}}
\end{eqnarray*}
satisfying the relation $\hat \psi_i(v_j)=\delta_{ij}$ for $i,j=1,\dots 4$. 
Corresponding gradients 
\begin{eqnarray*}
\nabla\hat \psi_1(x,y)&=&\left(-\frac{1}{h_x}+\frac{y}{h_{x}h_{y}}\,,\,-\frac{1}{h_y}+\frac{x}{h_{x}h_{y}}\right), \\
\nabla\hat \psi_2(x,y)&=&\left(\frac{1}{h_x}-\frac{y}{h_{x}h_{y}}\,,\,-\frac{x}{h_{x}h_{y}}\right), \\
\nabla\hat \psi_3(x,y)&=&\left(\frac{y}{h_{x}h_{y}}\,,\,\frac{x}{h_{x}h_{y}}\right), \\ 
\nabla\hat \psi_4(x,y)&=&\left(-\frac{y}{h_{x}h_{y}}\,,\,\frac{1}{h_y}-\frac{x}{h_{x}h_{y}}\right) 
\end{eqnarray*}
are linear functions.
Local stifness matrix 
is defined as
$$  ({\bf K}^{BIL}_{\rm{ref}})_{ij}=\int_{T_{\rm{ref}}} \nabla \hat \psi_i \cdot \nabla \hat \psi_j \,\textrm{d}x
$$
and direct computation shows 
\begin{eqnarray*}
{\bf K}^{BIL}_{\rm{ref}}= \frac{1}{6 h_x h_y}
\begin{pmatrix}
 2 h_x^2 + 2 h_y^2 & h_x^2 - 2 h_y^2 & - h_x^2 - h_y^2 & -2 h_x^2 + h_y^2 \\
 h_x^2 - 2 h_y^2 & 2 h_x^2 + 2 h_y^2 & -2 h_x^2 + h_y^2 & - h_x^2 - h_y^2 \\
- h_x^2 - h_y^2 & -2 h_x^2 + h_y^2  & 2 h_x^2 + 2 h_y^2 & h_x^2 - 2 h_y^2\\
-2 h_x^2 + h_y^2 & - h_x^2 - h_y^2 & h_x^2 - 2 h_y^2 & 2 h_x^2 + 2 h_y^2
\end{pmatrix}.
\end{eqnarray*}

Local Raviart-Thomas basis functions (of the lowest order) are vector edge based basis functions 
\begin{eqnarray*}
\hat \eta_1(x,y)&:=&\Bigl(0\,,\,1-\frac{y}{h_y}\Bigr), \\
\hat \eta_2(x,y)&:=&\Bigl(\frac{x}{h_x}\,,\,0\Bigr), \\
\hat \eta_3(x,y)&:=&\Bigl(0\,,\,\frac{y}{h_y}\Bigr), \\
\hat \eta_4(x,y)&:=&\Bigl(1-\frac{x}{h_x}\,,\,0\Bigr) 
\end{eqnarray*}
defined on reference edges 
$$e_1=\{v_1, v_2\}, \qquad e_2=\{v_2, v_3\}, \qquad e_3=\{v_3, v_4\}, \qquad e_4=\{v_4, v_1\}$$ and they satisfy the relation $\hat \eta_i |_{e_j} \cdot n_j = \delta_{ij},$
where global normals $n_j$ related to edges $e_j$ are always oriented in directions of the coordinate system, 
$$n_1=(0,1), \quad n_2=(1,0), \quad n_3=(0,1), \quad n_4=(1,0). $$ This choice of global normals leads to a simpler implementation with no issues related to global orientation of edges. The corresponding divergences 
\begin{eqnarray*}
{\rm{div}} \,\hat \eta_1 =-\frac{1}{h_y}, \quad 
{\rm{div}} \,\hat \eta_2=\frac{1}{h_x} , \quad 
{\rm{div}} \,\hat \eta_3=\frac{1}{h_y} , \quad 
{\rm{div}} \,\hat \eta_4=-\frac{1}{h_x} 
\end{eqnarray*}
are constant functions.
Local stifness and mass matrices defined by relations 
$$({\bf K}^{RT0 }_{\rm{ref}})_{ij}=\int_{T_{\rm{ref}}} {\rm{div}} \hat \eta_i  \; {\rm{div}} \hat \eta_j \,\textrm{d}x,   \qquad ({\bf M}^{\rm{RT0}}_{\rm{ref}})_{ij}=\int_{T_{\rm{ref}}} \hat \eta_i \cdot \hat \eta_j \,\textrm{d}x$$  
read
$${\bf K}^{RT0}_{\rm{ref}}=
\begin{pmatrix}
\frac{h_x}{h_y} & -1 & -\frac{h_x}{h_y} & 1 \\
-1 & \frac{h_y}{h_x} & 1 & -\frac{h_y}{h_x} \\
-\frac{h_x}{h_y} & 1 & \frac{h_x}{h_y} & -1 \\ 
1 & -\frac{h_y}{h_x} & -1 & \frac{h_y}{h_x}
\end{pmatrix}, \qquad {\bf M}^{RT0}_{\rm{ref}}=h_{x}h_{y}
\begin{pmatrix}
\frac{1}{3} & 0 & \frac{1}{6} & 0\\
0 & \frac{1}{3} & 0 & \frac{1}{6}\\
\frac{1}{6} & 0 & \frac{1}{3} & 0\\
0 & \frac{1}{6} & 0 & \frac{1}{3}
\end{pmatrix}.$$

\end{document}